\documentclass[journal]{new-aiaa}
\usepackage[utf8]{inputenc}
\usepackage{textcomp}

\usepackage{graphicx}
\usepackage{amsmath}
\usepackage[version=4]{mhchem}
\usepackage{siunitx}
\usepackage{longtable,tabularx}
\usepackage{tcolorbox}
\usepackage{caption}
\usepackage{float} 
\usepackage{subcaption}
\usepackage{threeparttable}
\usepackage{booktabs}
\usepackage{appendix}

\def\r{\overline{r}}

\newtheorem{lemma}{Lemma}

\newtheorem{problem}{Problem}

\renewcommand{\thesubsection}{\Alph{subsection}}

\title{Nonlinear Optimal Guidance for Intercepting Moving Targets}

\author{Han Wang \footnote{Ph.D. student, School of Aeronautics and Astronautics, Zhejiang.} and Zheng Chen\footnote{Researcher, School of Aeronautics and Astronautics, Zhejiang: z-chen@zju.edu.cn (Corresponding Author).}}
\affil{Zhejiang University, 310027 Hangzhou, People's Republic of China}
\affil{Huanjiang Laboratory, 311899 Zhuji, People's Republic of China}

\begin{document}

\maketitle

\begin{abstract}
This paper introduces a nonlinear optimal guidance framework for guiding a pursuer to intercept a moving target, with an emphasis on real-time generation of optimal feedback control for a nonlinear optimal control problem. Initially, considering the target moves without maneuvering, we derive the necessary optimality conditions using Pontryagin's Maximum Principle. These conditions reveal that each extremal trajectory is uniquely determined by two scalar parameters. Analyzing the geometric property of the parameterized extremal trajectories not only leads to an additional necessary condition but also allows to establish a sufficient condition for local optimality. This enables the generation of a dataset containing at least locally optimal trajectories. By studying the properties of the optimal feedback control, the size of the dataset is reduced significantly, allowing training a lightweight neural network to predict the optimal guidance command in real time. Furthermore, the performance of the neural network is enhanced by incorporating the target's acceleration, making it suitable for intercepting both uniformly moving and maneuvering targets. Finally, numerical simulations validate the proposed nonlinear optimal guidance framework, demonstrating its better performance over existing guidance laws.
\end{abstract}

\section{Introduction}\label{sec:Intro}
The Proportional Navigation (PN) is probably one of the most widely used guidance laws due to its simplicity and efficiency \cite{zarchan2012tactical}. It ensures that the pursuer’s acceleration is proportional to the Line-of-Sight (LOS)
rate, allowing for effective target interception even in scenarios of intercepting moving target \cite{4103772,102706}.  Recently, in order to satisfy various constraints, researchers have developed different variants of PN, such as biased PN with terminal angle constraint \cite{7849201} and biased PN for target observability enhancement \cite{7073471}. Furthermore, various methods have been used to develop advanced guidance laws against moving target \cite{9453157,doi:10.2514/1.G004762,8333799,10419004}. 

However, the methods mentioned above fail to consider optimality in terms of a meaningful performance index. By linearizing the kinematics around the collision triangle, a guidance law can be derived using linear-quadratic optimal control. Then, it is also used to find optimal PN. When the navigation constant is equal to 3, the PN has been mathematically proven to be optimal for intercepting nonmaneuvering targets in terms of control effort \cite{1975Applied}. However, when the deviations from collision triangle is relatively large, the control effort required by the PN may not be optimal; see, e.g., Ref.~\cite{chen2019nonlinear}. Thus, the methods in \cite{cho2014optimal} that are based on linear-quadratic optimal control inherently share the same limitations.

In recent decades, Nonlinear Optimal Guidance (NOG), by considering the nonlinear engagement kinematics, has been extensively studied, showing that it consumes less control effort than the PN in the nonlinear setting \cite{lu2006nonlinear}. However, the research on NOG is mainly focusing on problems with stationary targets; see, e.g., Refs.~\cite{lu2006nonlinear,chen2019nonlinear}. The current paper presents a natural extension to studying the NOG for intercepting moving targets. The fundamental problem for the NOG of intercepting a moving target is equivalent to finding the solution of a nonlinear optimal control problem within each guidance cycle or within a small period of time.

Up to now, the methods for solving nonlinear optimal control problems have been classified into two categories: indirect methods and direct methods \cite{2010Betts}. The indirect methods are based on Pontryagin's Maximum Principle (PMP) \cite{pontryagin1987mathematical}, which provides necessary conditions for optimality, and require solving two-point boundary value problems or multi-point boundary value problems \cite{1975Applied}. Although these kind of methods can produce precise solutions, they require iterative computations, making them impractical for real-time applications. 

On the other hand, in order to use numerical optimization techniques directly, the direct methods transform the optimal control problem to a parameter optimization problem \cite{2010joseph}. These years, various methods were devised, such as sequential convex programming-based optimal guidance \cite{kim2023optimal}, model predictive static programming based suboptimal guidance \cite{dwivedi2011suboptimal,10021226,10262163}, and geometric parameterization based optimal guidance \cite{10758229}. While these kinds of optimal guidance offer significant performance improvements, they often require intensive computational resources and are challenging to implement in real-time engagements. In order to realize the real-time generation of optimal trajectories, the convex programming based optimal guidance was developed; see, e.g., Refs.~\cite{10767281,8598772}. However, convexifying the nonlinear dynamics of intercepting a moving target can be extremely challenging.

Because of various issues of indirect methods and direct methods, the Neural Network (NN) has been combined in recent decades with optimal control methods to address optimal control problems \cite{doi:10.2514/1.G005254,9805687,9163239,8587201,3551283} in a real-time manner. By training an NN on a dataset of optimal trajectories obtained from solving the nonlinear optimal control problem offline, we can approximate the nonlinear optimal guidance law and achieve real-time performance. However, the dataset generate by indirect and direct methods cannot be guaranteed to be optimal \cite{doi:10.1137/15M1013274,doi:10.2514/1.G005865}. 

To address the above issues, Ref.~\cite{doi:10.2514/1.G006666} proposed a parameterized system for generating optimal trajectories using PMP in addition to some extra optimal conditions. By simply solving some initial value problems, Wang et al \cite{doi:10.2514/1.G006666} generated a set of solutions to the optimal control problem of intercepting stationary target. This approach allows for the creation of a dataset mapping the pursuer's state to the corresponding NOG command. As a continuation of \cite{doi:10.2514/1.G006666}, this paper extends to study a more important and complex nonlinear optimal guidance problem for intercepting moving targets, ensuring both computational efficiency and optimal control performance. According to the PMP, it is found in the paper that the extremal trajectories are determined by two scalar parameters. Then, by embedding sufficient conditions into the parameterized extremal trajectories, one is able to construct the dataset of at least locally optimal trajectories for the nonlinear optimal guidance problem with moving targets. Training a lightweight neural network by the dataset eventually allows to generate at least locally nonlinear optimal guidance command for intercepting moving targets. 

By embedding the above method into the closed-loop guidance system, it can be extended to generating NOG for intercepting maneuvering targets. Many existing augmented guidance methods, such as Augmented Proportional Navigation (APN) \cite{zarchan2012tactical}, Sliding Mode Guidance (SMG) \cite{7472963}, and Pseudocontrol-Effort Optimal Guidance (PEOG) \cite{jeon2015exact}, have been developed to enhance performance against maneuvering targets. By incorporating augmentation terms inspired by these methods, the trained neural network can be further refined to improve robustness against unpredictable target maneuvers. This integration enables the guidance law to dynamically adapt to different target behaviors while maintaining near-optimal control performance.

The remaining sections of this paper are structured as follows: Section~\ref{Sec:PF} formulates the nonlinear optimal control problem for intercepting a nonmaneuvering target. In Section~\ref{Sec:OGNN}, after deriving the necessary conditions of optimality from PMP, a parameterized set of differential equations is introduced to describe the optimal solutions. Additional necessary and sufficient conditions for local optimality are presented. Section~\ref{Sec:GSD} discusses the optimal guidance architecture and the geometry properties of the control command. The procedure of generating the sampled dataset for training the neural network is also detailed. In Section~\ref{sec:GMT}, a guidance law for intercepting maneuvering targets is proposed, incorporating an augmentation term. Finally, in Section~\ref{Sec:NS}, numerical simulations are presented to illustrate the effectiveness of the proposed approach.

%To address the limitations of PN, researchers have developed advanced guidance laws based on nonlinear control theories. From the perspective of nonlinear control, a guidance problem can be regard as a finite-time tracking problem. In order to intercept the target, zero-effort-miss \cite{doi:10.2514/1.14951,doi:10.2514/1.24953,7812864,doi:10.2514/1.G002949} or LOS rate \cite{doi:10.2514/2.4792,53460,doi:10.2514/1.42976,7472963} is devised as the tracking error to be nullified. These methods offer significant advantages, such as strong disturbance rejection, adaptability to changing engagement conditions, and improved handling of nonlinear target dynamics. However, they also come with drawbacks, including lack of explicit optimality criteria \cite{doi:10.2514/1.G003343} and potential chattering effects in sliding mode control \cite{doi:10.2514/2.4792}. In recent years, researchers have attempted to consider the optimality of problems when using the nonlinear control theories. For example, \cite{doi:10.2514/1.G003343} derived an optimal error dynamics for a general form of guidance problem. However it is hard to design the error for a nonlinear guidance problem to tracking the error dynamics. Therefore, while these methods can improve robustness, they still do not explicitly minimize control effort, leading to inefficient energy consumption. 
\section{Problem Formulation}\label{Sec:PF}
\begin{figure}[bt!]
	\centering	
	\includegraphics[width=.4\textwidth]{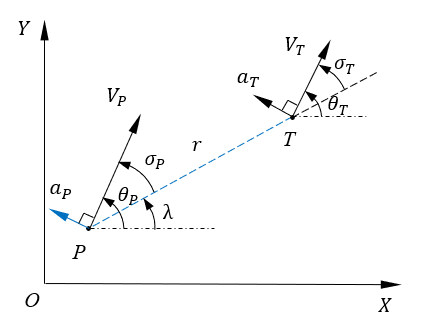}
	\caption{Engagement geometry.}
	\label{fig:Engage}
\end{figure}

Let us consider a two-dimensional interception problem in an inertial Cartesian coordinate frame $OXY$, as shown in Fig.~\ref{fig:Engage}, and we denote by P and T the pursuer and the target, respectively. The speeds of the pursuer and the target are assumed to be constant, and are represented by $V_P$ and $V_T$, respectively. The heading angles for the pursuer and the target are given by $\theta_P\in [0,2\pi)$ and $\theta_T\in [0,2\pi)$, respectively, which are positive when measured counterclockwise. Let $\lambda\in [0,2\pi)$ be the angle between $OX$ axis and the LOS.
Denote by $\sigma_P\in [-\pi,\pi)$ and $\sigma_T\in [-\pi,\pi)$ the lead angles of the pursuer and the target, respectively, and they are expressed as
\begin{align}
	\sigma_P &= \theta_P - \lambda\nonumber\\
	\sigma_T &= \theta_T - \lambda
\end{align}
Let $\r > 0$ be the Euclidean distance between the pursuer and the target. Then, the differential equations governing the relative motions of the pursuer and the target can be expressed as  %The positions of the pursuer and the target in $X_IOY_I$ are given by $(x_P,y_P)\in \mathbb{R}^2$ and $(x_T,y_T)\in \mathbb{R}^2$, respectively. 
\begin{equation}
	\label{eq:r-lambda}
	\begin{split}
		&\dot \r=-V_P \cos\sigma_P+V_T \cos\sigma_T\\
		&\dot \lambda=-\frac{V_P}{\r} \sin\sigma_P+\frac{V_T}{\r} \sin\sigma_T\\
		&\dot \theta_P=\frac{a_p}{V_P}\\
		&\dot \theta_T=\frac{a_T}{V_P}
	\end{split}
\end{equation}
where the over dot denotes the differentiation with respect to time, $a_P\in \mathbb{R}$ and $a_T\in \mathbb{R}$ are the lateral accelerations of the pursuer and the target, respectively.

By normalizing the speed of pursuer to one and assuming that the target moves without maneuvering, the kinematics can be simplified to
\begin{equation}
	\label{eq:norm-r-lambda}
	\begin{cases}
		\dot r(t)=-\cos(\theta_P-\lambda)+\mu\cos(\theta_T-\lambda)\\
		\dot \lambda(t)=\displaystyle{\frac{-\sin(\theta_P-\lambda)+\mu\sin(\theta_T-\lambda)}{r}}\\
		\dot \theta_P(t)=u\\
		%&\dot \theta_T=0
	\end{cases}
\end{equation}
where $r$ is normalized distance between the pursuer and the target, the constant $\mu=V_T/V_P$ is the speed ratio, and $u\in \mathbb{R}$ is the control parameter, related to the lateral acceleration of the pursuer. Note that if $\mu = 0$, it is related to the problem of intercepting a stationary target, which has been studied in
\cite{chen2019nonlinear}. It should also be noted that the scenario of $\mu\geq1$ is usually not considered for the intercepting problem; see, e.g., \cite{ghosh2014capturability,jeon2015exact,doi:10.2514/1.G005868}. Thus, we assume $\mu \in (0,1)$ in the remainder of this paper. %Because the maneuver of target is not considered, the normalized relative normal acceleration $u=a_P/V_P$. The LOS angle in frame $X_ROY_R$ is denoted as $\bar\lambda$.

According to definitions and notations above, the NOG problem for intercepting a moving  but nonmaneuvering target is equivalent to the following Optimal Control Problem (OCP). 
\begin{problem}[OCP]\label{pro1}
	Given a pursuer and a moving but nonmaneuvering target, let the initial heading angles of the pursuer and the target be $\theta_{P0}$ and $\theta_{T0}$, respectively,  let the initial distance between the pursuer and the target be $r_0$,  the initial LOS angle be $\lambda_0$, and let the speed ratio $\mu$ take a value in $(0,1)$. Then,  the OCP consists of finding a measurable control $(0,t_f)\mapsto u$ that steers the system in Eq. (\ref{eq:norm-r-lambda}) from the initial condition $(r_0,\lambda_0,\theta_{P0})$ to intercepting the moving target, i.e., $r(t_f) = 0$, so that
	\begin{equation*}
		%\label{sample:equation7}
		J = \int_{0}^{t_f} \left[\kappa+(1-\kappa)\frac{1}{2}u(t)^2\right]dt
	\end{equation*}
	is minimized where $t_f$ is the free final time and $\kappa\in(0,1)$ is a weighting factor.
\end{problem}
%In order to make sure the solution of the OCP exists, it is required that the speed ratio $\mu$ takes values in $[0,1)$, as analyzed in \cite{?}. 
% By the following lemma, we shall show that the solution of the OCP exists if the speed ratio $\mu$ takes values in $(0,1)$.
%\begin{lemma}\label{lemma:exist}
%	Given any initial condition $(r_0,\lambda_0,\theta_{P0})$ and any weighting factor $\kappa \in (0,1)$, if the moving target is not maneuvering with $\mu \in (0,1)$, then the OCP has a global optimum. 
%\end{lemma}
%The proof is postponed to Appendix \ref{app:0}.

During guiding a pursuer to intercepting a moving target, it is strictly required that the onboard computer of the pursuer produces the optimal guidance command within each guidance cycle. To this end, the OCP in Problem \ref{pro1} should be solved in real time or within a small period of time. However, as stated in the Introduction, existing numerical methods cannot guarantee to find the solution of the OCP in real time. In the subsequent sections, a parameterized approach will be presented for obtaining the optimal guidance command in real time via neural network. 

\section{Characterization of Optimal Trajectories}\label{Sec:OGNN}
In this section, we first present some necessary conditions for optimality from PMP and then use these necessary conditions to establish a parameterized family of extremals.% Subsequently, we elucidate the geometric properties inherent to the optimal control. Furthermore, leveraging these conditions, we parameterize the optimal solutions using two bounded scalars. For brevity, key results and propositions are encapsulated within lemmas, with their proofs deferred to the appendices.%
\subsection{Necessary Conditions}
Denote by $(p_r,p_\lambda,p_{\theta})\in\mathbb{R}^3$ the costate of the state $(r,\lambda,\theta_P)$. Then, the Hamiltonian for the OCP is expressed as
\begin{equation}
	\label{sample:Ham}
	\begin{split}
		H = p_r(-\cos\sigma_P+\mu\cos\sigma_T)+p_\lambda\frac{-\sin\sigma_P}{r}+p_\lambda\frac{\mu\sin\sigma_T}{r}+p_{\theta}u+p^0\left[\kappa+\frac{1}{2}(1-\kappa)u^2\right]
	\end{split}	
\end{equation}
where $p^0$ is a negative scalar according to \cite[Remark 2]{chen2019nonlinear}. Because for any negative $p^0$ the quadruple $(p_x,p_y,p_\theta,p^0)$ can be normalized so that $p^0 = -1$, we shall consider $p^0 = -1$ in the remainder of the paper.

The costate variables are governed by 
\begin{equation}
	\label{eq:ode1}
	\begin{cases}
		\dot p_r(t) = p_\lambda\frac{-\sin\sigma_P+\mu\sin\sigma_T}{r^2}\\
		\dot p_\lambda(t) = p_r(\sin\sigma_P-\mu\sin\sigma_T)+p_\lambda\frac{\cos\sigma_P-\mu\cos\sigma_T}{r}\\
		\dot p_\theta(t) = -p_r\sin\sigma_P+p_\lambda\frac{\cos\sigma_P}{r}
	\end{cases}
\end{equation}
According to PMP \cite{pontryagin1987mathematical}, we have 
\begin{equation}
	\label{eq:u}
	\frac{\partial H}{\partial u}=0
\end{equation}
%\emph{Remark 2}: When $p^0 = 0$, the explicit formula of (\ref{eq:u}) implies $p_\theta \equiv 0$, which indicates $\dot {p_\theta} \equiv 0$. According to the third equation of (\ref{eq:ode1}), if $\dot p_\theta \equiv 0$, the optimal trajectory of the OIP is a straight line, and hence the corresponding optimal control is null, which happens only if the initial velocity vector points to the target (or the initial heading angle $\theta_0$ is $\pi$). Because $p^0$ is nonpositive, and because we will not consider the interval of $\sigma_M$ to be $0$ (see Remark 1), we have that $p^0$ is negative. For any negative $p^0$, the quadruple $(p_x,p_y,p_\theta,p^0)$ can be normalized so that $p^0 = -1$. Thus, we shall consider $p^0 = -1$ in the remainder of the paper.
which can be written explicitly as
\begin{equation}
	%\label{eq:u}
	u(t)=\frac{p_\theta}{1-\kappa}
\end{equation}
Because $\lambda(t_f)$ and $\theta_P(t_f)$ are not fixed, the transversality condition implies
\begin{equation}
	\label{equation:transversality}
	\begin{cases}
		p_\lambda(t_f)= 0\\
		p_\theta(t_f)= 0
	\end{cases}	
\end{equation}
As the final time is free, we have 
\begin{equation}
	\label{equation:Hamiltonian}
	H \equiv 0
\end{equation}
along any optimal trajectory. 

For notational simplicity, we refer to a triple $(r(t),\lambda(t),\theta_P(t))$ for $t\in [0,t_f]$ as an extremal trajectory if it satisfies all the necessary conditions given in Eqs. (\ref{sample:Ham}-\ref{equation:Hamiltonian}). In addition, the control along an extremal trajectory will be said as extremal control. In the following subsection, we shall establish a parameterized family of extremal trajectories.

%%%%%%%%%%%%%%%%%%%%%%%%%%%%%%%%%%%%%%%%%%%%%%%%%%%%%
%\section{Parametrization of Extremal Trajectories}\label{Sec:SOC}
%%%%%%%%%%%%%%%%%%%%%%%%%%%%%%%%%%%%%%%%%%%%%%%%%%%%%
%In this section, we first parameterize the aforementioned necessary conditions. Then, to decrease the one-to-many mapping, we establish a parametrization approach for the Hamiltonian extremals satisfying the sufficient conditions. Furthermore, in order to detecting conjugate points, the method for obtaining the value of $\delta(\bar t)$ is presented.

\subsection{Parametrization of Extremal Trajectories}\label{sec:HE}

%Let $R$ be in $\mathbb{R}_0^+$, and let $\Lambda$ and $\Theta$ take values in $[0,2\pi)$. Let $(P_R,P_\Lambda,P_\Theta)\in\mathbb{R}^3$. 
Given $R\geq 0$, $\Lambda \in [0,2\pi)$,  $\Theta \in [0,2\pi)$, and $(P_R,P_\Lambda,P_\Theta)\in \mathbb{R}^3$, let us introduce the following differential equations:
%\begin{figure*}[htp]
	\begin{equation}
		\label{equation:phpx1}
		\begin{cases}
			\dot R(\tau)=\cos(\Theta-\Lambda)-\mu\cos\Lambda\\
			\dot \Lambda(\tau)=\displaystyle{\frac{\sin(\Theta-\Lambda)+\mu\sin\Lambda}{R}}\\
			\dot \Theta(\tau)=-\displaystyle{\frac{P_\Theta}{1-\kappa}}\\
			\dot P_R(\tau) = -P_\Lambda\displaystyle{\frac{-\sin(\Theta-\Lambda)-\mu\sin\Lambda}{R^2}}\\
			\dot P_\Lambda(\tau) = -P_R(\sin(\Theta-\Lambda)+\mu\sin\Lambda)-P_\Lambda\displaystyle{\frac{\cos(\Theta-\Lambda)-\mu\cos\Lambda}{R}}\\
			\dot P_\Theta(\tau) = P_R\sin(\Theta-\Lambda)-P_\Lambda\displaystyle{\frac{\cos(\Theta-\Lambda)}{R}}
		\end{cases}
	\end{equation}
%\end{figure*}
where the initial values at $\tau=0$ are set to satisfy the following equations
\begin{equation}
	\begin{cases}
		R(0) = 0\\
		P_{\Lambda}(0) = 0\\
		P_{\Theta}(0) = 0\\
		P_{R}(0) [-\cos (\Theta(0) - \Lambda(0)) + \mu \cos \Lambda(0)]= \kappa
	\end{cases}
	\label{eq:H=0}
\end{equation}
Because $\kappa\neq0$, Eq. (\ref{eq:H=0}) indicates
\begin{equation}
	\cos(\Theta(0)-\Lambda(0))\neq\mu\cos\Lambda(0)
\end{equation}
Thus, by solving Eq. (\ref{eq:H=0}), $P_{R}(0)$ can be expressed as a function of $\Theta(0)$ and $\Lambda(0)$, i.e., 
\begin{equation}
	P_{R}(0) = \frac{\kappa}{-\cos(\Theta(0)-\Lambda(0))+\mu\cos\Lambda(0)}
\end{equation}
Up to now, it has been apparent that, given any speed ratio $\mu\in (0,1)$, the solution of Eq. (\ref{equation:phpx1}) with the initial condition given in Eq. (\ref{eq:H=0}) at any $\tau \geq 0$ is totally determined by $\Lambda(0)$ and $\Theta(0)$. Thus, if denoting  $\Lambda(0)$ and $\Theta(0)$ by $\Lambda_0$ and $\Theta_0$, respectively, we have that for any given $\mu \in (0,1)$ the solution of the initial value problem defined in Eq. (\ref{equation:phpx1}) and Eq. (\ref{eq:H=0}) is totally determined by the parameters $\tau$, $\Lambda_0$, and $\Theta_0$. For notational simplicity, given any speed ratio $\mu \in (0,1)$ we denote by 
\begin{equation*}
	\begin{split}
		\mathcal{T}^\mu(\tau,\Theta_0,\Lambda_0):=(R^\mu(\tau,\Theta_0,\Lambda_0),\Lambda^\mu(\tau,\Theta_0,\Lambda_0),\Theta^\mu(\tau,\Theta_0,\Lambda_0),P^\mu_R(\tau,\Theta_0,\Lambda_0),P^\mu_\Lambda(\tau,\Theta_0,\Lambda_0),P^\mu_\Theta(\tau,\Theta_0,\Lambda_0))
	\end{split}
\end{equation*}
the solution of the initial value problem in Eq.~(\ref{equation:phpx1}) and Eq.~(\ref{eq:H=0}). In addition, we use $\Pi$ to denote the projection from the cotangent bundle to the state space, i.e.,
\begin{equation*}%\label{eq:Pi}
	\begin{split}
		\Pi(\mathcal{T}^\mu(\tau,\Theta_0,\Lambda_0)) = (R^\mu(\tau,\Theta_0,\Lambda_0),\Lambda^\mu(\tau,\Theta_0,\Lambda_0),\Theta^\mu(\tau,\Theta_0,\Lambda_0))
	\end{split}	
\end{equation*} 

\begin{lemma}\label{lemma:ex_tra}
	Given a pursue and a target, let the heading angle of the target be 0, i.e., $\theta_T = 0$. Then, for any speed ratio $\mu \in (0,1)$ and any initial condition $(r_0,\lambda_0,\theta_{P0})$ for Problem \ref{pro1}, there exists $\tau>0$, $\Theta_0 \in [0,2\pi)$, and $\Lambda_0\in [0,2\pi)$ so that
		$$\Pi(\mathcal{T}^{\mu}(\tau,\Theta_0,\Lambda_0) ) = (r_0,\lambda_0,\theta_{P0})$$
	Conversely, given any $\mu \in (0,1)$ and any $\Theta_0,\Lambda_0 \in [0,2\pi)$, there exists an initial condition $(r_0,\lambda_0,\theta_{P0})$ for Problem \ref{pro1} so that the solution trajectory $\Pi(\mathcal{T}^{\mu}(\cdot,\Theta_0,\Lambda_0))$ on $[0,t_f]$ is the reverse of the optimal trajectory of Problem \ref{pro1}, i.e.,
		$$(r(t),\lambda(t),\theta_P(t)) =   \Pi(\mathcal{T}^{\mu}(t_f-t ,\Theta_0,\Lambda_0))  $$
	where $(r(t),\lambda(t),\theta_P(t))$ is the optimal trajectory of Problem \ref{pro1}. 
\end{lemma} 
The proof is postponed to Appendix \ref{app:A}.

Set
\begin{equation}
	U^\mu(\tau,\Theta_0,\Lambda_0):=\frac{P^\mu_\Theta(\tau,\Theta_0,\Lambda_0)}{1-\kappa}
\end{equation}
%Denote by $u^*(x,y,\theta)$ the optimal control at the state $(x,y,\theta)$. 
It is evident that $U^\mu(\tau,\Theta_0,\Lambda_0)$ represents the extremal control along the  extremal trajectory $\Pi(\mathcal{T}^\mu(\tau,\Theta_0,\Lambda_0))$. %the extremal trajectory
As a result of Lemma \ref{lemma:ex_tra}, one can use the initial value problem defined by Eq. (\ref{equation:phpx1}) and Eq. (\ref{eq:H=0}) to generate the dataset of extremal trajectories for the specific case that the heading angle of the target is zero. This will be vital in the next section to establish the closed-loop optimal guidance scheme. Before proceeding to the next section, we supplement some additional optimality conditions by analyzing the properties of the parameterized extremals in the following subsection. 

\subsection{Supplementary Optimality Conditions}
%It is evident that the necessary optimality conditions from the PMP are incorporated in the initial value problem defined by (\ref{equation:phpx1}) and (\ref{eq:H=0}). In this subsection, we shall present another necessary condition for optimality by analyzing the geometric property of trajectories and a sufficient condition for local optimality. 

By analyzing the geometry property of extremal trajectory, we present an extra optimality condition by the following lemma. 

\begin{lemma}
	\label{lemma:3}
	Given any trajectory $\Pi(\mathcal{T}^\mu(\tau,\Theta_0,\Lambda_0))$ on $[0,\tau_f]$ with $\tau_f$ being a positive number, if there exists a time $\bar\tau$ within interval $(0,\tau_f)$ such that the velocity vector $[\cos\Theta^\mu(\bar\tau,\Theta_0,\Lambda_0),\sin\Theta^\mu(\bar\tau,\Theta_0,\Lambda_0)]$ is collinear with the LOS, i.e., 
	\begin{equation}\label{eq:sin_dth}
		\Theta^\mu(\bar\tau,\Theta_0,\Lambda_0)=\Lambda^\mu(\bar\tau,\Theta_0,\Lambda_0)
	\end{equation}
	then the trajectory $\Pi(\mathcal{T}^\mu(\tau,\Theta_0,\Lambda_0))$ on $[0,\tau_f]$ is not optimal.
\end{lemma}
The proof is postponed to Appendix \ref{app:B}. 

It is apparent that Lemma \ref{lemma:3} supplements an additional necessary condition for optimality. It is important to note that these necessary conditions alone do not ensure that a solution trajectory is at least locally optimal unless additional sufficient conditions hold. The following lemma provides a sufficient condition for establishing local optimality.

\begin{lemma}\label{lemma:2}%[Sufficient Optimality Conditions]
	Given any trajectory $\Pi(\mathcal{T}^\mu(\tau,\Theta_0,\Lambda_0))$ on $[0,\tau_f]$ with $\tau_f$ being a positive number, set
	\begin{equation}\label{eq:lemma1}
		\delta(\tau,\Theta_0,\Lambda_0) :=\mathrm{det}\left[\frac{\partial\Pi(\mathcal{T}^\mu(\tau,\Theta_0,\Lambda_0))}{\partial(\tau,\Theta_0,\Lambda_0)}\right]
	\end{equation}
If $\delta(\tau,\Theta_0,\Lambda_0)\neq0$ on $\tau\in(0,\tau_f]$, then the trajectory $\Pi(\mathcal{T}^\mu(\tau,\Theta_0,\Lambda_0))$ on $[0,\tau_f]$ is a local optimum; if the determinant $\delta(\tau,\Theta_0,\Lambda_0)$ changes its sign at a time $\bar{\tau}\in (0,\tau_f)$, then the trajectory $\Pi(\mathcal{T}^\mu(\tau,\Theta_0,\Lambda_0))$ for $\tau \in [0,\tau_f]$ loses its local optimum.
\end{lemma}
Lemma \ref{lemma:2} presents a sufficient condition for local optimality, and it is a direct result from Theorem 1 in \cite{doi:10.2514/1.G000284,doi:10.2514/1.G005865}. Readers who are interested in the corresponding proof are referred to \cite{doi:10.2514/1.G000284,doi:10.2514/1.G005865}.

Up to now, it has been evident that given any $\Theta_0$ and $\Lambda_0$, the extremal trajectory $\Pi(\mathcal{T}^\mu(\tau,\Theta_0,\Lambda_0))$ for $\tau \in (0,\tau_f)$ satisfying the necessary condition in Lemma \ref{lemma:3} and and the sufficient condition, i.e., $\delta(\tau,\Theta_0,\Lambda_0)\neq 0$ for $\tau\in (0,\tau_f]$ in Lemma \ref{lemma:2}, is at least a locally optimal solution. In the next section, these optimality conditions will be employed to generate the optimal guidance command for intercepting moving but nonmaneuvering targets.
%Because the parameterized system $\mathcal{F}$ only satisfy the necessary conditions from PMP, the dataset generated from $\mathcal{F}$ will contain incorrect control, which means one-to-many mapping (see Remark 1). To overcome this drawback, the sufficient conditions are established by leveraging the conjugate point theorem from the geometry optimal control theory. With these improved conditions of optimality, we are able to generate the dataset that is at least locally optimal.
\section{NOG for Intercepting Nonmaneuvering Targets}\label{Sec:GSD}
%Due to the complexity of trajectories of aerospace vehicles, the DNN was usually used for the representation of the relationship between the state and the corresponding optimal control. However, it is a quite burden for the computer onboard because of the large scale of the neural network. Thus, in the following paragraphs, we shall show how to generate the feedback optimal control with a simple feedforward neural network. By the kinematics in (\ref{eq:ode1}) and the symmetric property developed in last section, the space of optimal solutions is shrunk.
In this section, we first establish the guidance architecture of the pursuer. Then, we shall show how to generate the optimal guidance command with a simple neural network using the aforementioned parameterized system.

\subsection{Guidance Architecture}
Note that the optimal feedback control is not only determined by the state $(r,\lambda,\theta_P)$ of the pursuer, but also affected by the speed ratio $\mu \in (0,1)$ and the heading angle $\theta_T$ of the target. Thus, let us denote by $u^*(r,\lambda,\theta_P,\theta_T,\mu)$ the optimal feedback control at the state $(r,\lambda,\theta_P)$ for the pursuer to intercept the target with its heading angle being $\theta_T$. Then, given an optimal trajectory $(r(t),\lambda(t),\theta_P(t))$ of the OCP for a specified $\theta_T$ and $\mu$, let $u(t)$ denote the corresponding time history of optimal control. Then, for any $t\in [0,t_f]$ the following equation holds:
\begin{align}
	u^*(r(t),\lambda(t),\theta_{P}(t),\theta_T,\mu)=u(t)
\end{align}
Consequently, addressing the OCP in Problem \ref{pro1} in real time is equivalent to finding the value of the optimal feedback control $u^*(r,\lambda,\theta_P,\theta_T,\mu)$ for any $(r,\lambda,\theta_P,\theta_T,\mu)$ in real time. 

%Let $t_c\in[0,t_f)$ be the current time and $(r_c,\lambda_c,\theta_{Pc})$ be the state at $t_c$. Considering different OCPs with different $\theta_T\in[0,2\pi)$ and $\mu\in(0,1)$, the optimal control is not only the function of the current state $(r_c,\lambda_c,\theta_{Pc})$ but also related with the heading angle of the target $\theta_T$ and the speed ratio $\mu$. Thus, given an optimal trajectory $(r(t),\lambda(t),\theta_P(t))$ of the OCP with heading angle of the target $\theta_T$ and the speed ratio $\mu$, if $u(t)$ is the corresponding optimal control, then for any $t\in[0,t_f)$, the following equation holds

%In result, addressing the OCP is equivalent to obtaining the value of the optimal feedback control $u^*(r_c,\lambda_c,\theta_{Pc},\theta_T,\mu)$ for any $(r_c,\lambda_c,\theta_{Pc},\theta_T,\mu)$ in real time.

\begin{figure}[hbt!]
	\centering	
	\includegraphics[width=.6\textwidth]{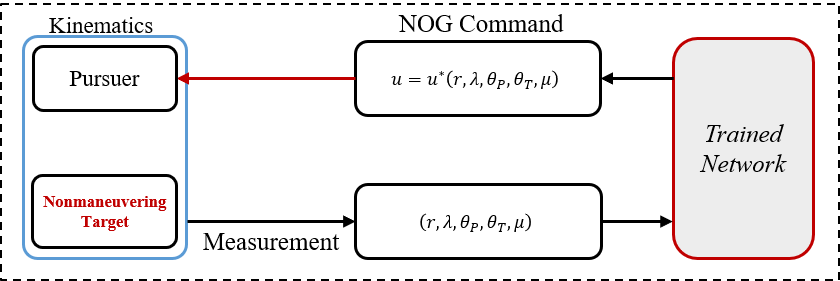}
	\caption{Guidance framework for intercepting a nonmaneuvering target.}
	\label{fig:Guidance}
\end{figure}

According to the universal approximation theorem \cite{1989Multilayer}, if a dataset capturing the mapping $(r,\lambda,\theta_P,\theta_T,\mu)\mapsto u^*(r,\lambda,\theta_P,\theta_T,\mu)$ can be obtained, an NN can be trained by the dataset to approximate it. Due to the simple structure of NN, its output is simply a composition of multiple linear transformations applied to the input vector. Thus, given an input $(r,\lambda,\theta_P,\theta_T,\mu)$, the output of the trained NN that represents the optimal feedback control $u^*(r,\lambda,\theta_P,\theta_T,\mu)$ can be obtained within a constant time. If the trained NN is embedded in the closed-loop guidance system, as shown in Fig.~\ref{fig:Guidance}, it will play the role of generating the optimal guidance command for intercepting nonmaneuvering target. According to the above analysis, the core of using an NN for generation of optimal guidance command lies in first constructing the dataset that captures the desired  mapping $(r,\lambda,\theta_P,\theta_T,\mu)\mapsto u^*(r,\lambda,\theta_P,\theta_T,\mu)$. In the following subsection, the procedure for constructing the dataset will be presented by applying the developments in Section \ref{Sec:OGNN}. 

\subsection{Dataset Generation for Training NN}
Thanks to Lemma \ref{lemma:ex_tra}, we can generate the dataset for the mapping $(r,\lambda,\theta_P,0,\mu)\mapsto u^*(r,\lambda,\theta_P,0,\mu)$ by solving the initial value problem defined in Eq. (\ref{equation:phpx1}) and Eq. (\ref{eq:H=0}) through sampling the values of $\Theta_0$, $\Lambda_0$, and $\mu$. Given any $\mu \in (0,1)$, an initial condition $(r_0,\lambda_0,\theta_{P0})$ is said to be feasible for Problem \ref{pro1} if there exists a trajectory starting from $(r_0,\lambda_0,\theta_{P0})$ to the moving target. Then, the rotation property in the following lemma shall show that it is enough to fix $\theta_T$ as zero during constructing the dataset. 

\begin{lemma}\label{lemma:rotation}
	Given any feasible state $(r,\lambda,\theta_{P})$ of the pursuer, a heading angle $\theta_{T}\in(0,2\pi)$ of the target, and a speed ratio $\mu\in(0,1)$, we have 
	\begin{equation}\label{eq:rotation}
		u^*(r,\lambda,\theta_{P},\theta_{T},\mu)=u^*(r,\lambda+\theta_{T},\theta_{P}+\theta_{T},0,\mu)
	\end{equation}
\end{lemma}
The proof is postponed to Appendix \ref{app:C}. 

Let us denote by $C(r,\lambda,\theta_{P},\mu)$ the mapping from $(r,\lambda,\theta_{P},0,\mu)$ to $u^*(r,\lambda,\theta_{P},0,\mu)$. Then, Lemma \ref{lemma:rotation} indicates that for any $(r,\lambda,\theta_{P},\theta_{T},\mu)$ with $\theta_{T}\neq0$, we have
\begin{equation*}
	u^*(r,\lambda,\theta_{P},\theta_{T},\mu)=C(r,\lambda+\theta_T,\theta_{P}+\theta_T,\mu)
\end{equation*}

\begin{lemma}\label{lemma5}
	Let the heading angle of the moving target be zero, i.e., $\theta_T = 0$. Then, given any $\mu \in (0,1)$ and any positive number $R_B > 0$, we have that for any feasible state $(r,\lambda,\theta_P)$, there exists $\tau>0$, $(\Theta_0,\Lambda_0)\in [0,2\pi)^2$, and a positive scalar $s> 0$ so that 
	\begin{equation*}
		\begin{cases}
			\displaystyle{\frac{r}{s}} = R^\mu(\tau,\Theta_0,\Lambda_0)\\
			\lambda = \Lambda^\mu(\tau,\Theta_0,\Lambda_0)\\
			\theta_P = \Theta^\mu(\tau,\Theta_0,\Lambda_0)\\
			C(r,\lambda,\theta_P,\mu) = \displaystyle{\frac{1}{s}}C(\frac{r}{s},\lambda,\theta_P,\mu)\\
		\end{cases}
	\end{equation*}
\end{lemma}
The proof is postponed to Appendix \ref{app:D}. 

Denote by $\mathcal{F}^{\mu}$ the set of all the feasible state. Then, for any $\mu \in (0,1)$, let us define a subset of feasible set as 
$$\mathcal{F}^{\mu}(R_B):= \{(r,\lambda,\theta_P)\in \mathcal{F}^{\mu} | r < R_B\}$$
where $R_B$ is a positive number. Then, Lemma \ref{lemma5} indicates that for any feasible state $(r,\lambda,\theta_P)\not\in \mathcal{F}^{\mu}(R_B)$ with $\theta_T$ and $\mu$, we can find a scalar $s>1$ so that $(r/s,\lambda,\theta_P)\in \mathcal{F}^{\mu}(R_B)$ and
$$C(r,\lambda+\theta_T,\theta_P+\theta_T,\mu) =\displaystyle{\frac{1}{s}}C(\frac{r}{s},\lambda+\theta_T,\theta_P+\theta_T,\mu)$$

Notice that one is able to gather the dataset by solving the initial value problem defined in Eq. (\ref{equation:phpx1}) and Eq. (\ref{eq:H=0}) via sampling  $(\Theta_0,\Lambda_0)\in [0,2\pi)^2$. However, the properties established in Lemma \ref{lemma:rotation} and Lemma \ref{lemma5} allows to significantly reduce the size of dataset. Lemma \ref{lemma:rotation} indicates that it is enough to set $\theta_T$ as zero during constructing the dataset, and Lemma \ref{lemma5} indicates that it is enough to choose the dataset with $r$ smaller than a positive number $R_B$. The detailed procedure for constructing the dataset is summarized in {\it Procedure \ref{pro1}}. 

\begin{center}
	\begin{tcolorbox}[
		colframe=blue!25,
		colback=gray!20,
		coltitle=blue!20!black,  
		fonttitle=\bfseries,
		width = 1\columnwidth,
		boxrule=0pt,
		top=1pt,
		bottom=-0.5pt,
		%enhanced,
		%segmentation style={solid, black,line width=1.2pt},
		segmentation code={\draw[black,solid,line width=1.2pt]($(segmentation.west)+(0.2,0)$)--($(segmentation.east)+(-0.2,0)$);}]
		\begin{center}
			\textbf{Procedure 1}: {\it Generation of Sampled Data for Optimal Feedback Control}
		\end{center}
		%\tcblower
		\bigskip
		
		\begin{itemize}
			\item[1.] Let $(\Theta_0^i,\Lambda_0^i)$, $i=1,2,\ldots,N$, be the uniformly chosen points from $[0,2\pi)^2$. Let $\mu^j$, $j=1,2,\ldots,M,$ be the uniformly chosen values of speed ratio from $(0,1)$. Let $R_B$ and $h$ be two positive numbers. Set $i=1$, $j=1$, and $\mathcal{D}= \varnothing$. 
			\item[2.] If $i\leq N$, go to step 3; otherwise, go to step 7.
			\item[3.] If $j\leq M$, set $\bar t=0$ and go to step 4; otherwise, set $i=i+1$ and go to step 2.
			%\item[4.] If $\bar t+h \leq T$, go to step 5; otherwise, set $j=j+1$ and go to step 3.	
			\item[4.] Propagate the system in Eq. (\ref{equation:phpx1}) with $(\Theta_0^i,\Lambda_0^i)$ and $\mu^j$ to generate the extremal $\Pi(\mathcal{T}^\mu(t,\Theta_0^i,\Lambda_0^i))$ and extremal control $U^\mu(t,\Theta_0^i,\Lambda_0^i)$ for $t\in[\bar t,\bar t+h]$, and go to step 5.
			\item[5.] If the determinant $\delta(\bar t,\Theta_0^i,\Lambda_0^i)$ in Lemma \ref{lemma:2} changes its sign or the supplemented necessary condition in Lemma \ref{lemma:3} is not satisfied, or $R^\mu(\bar t,\Theta_0^i,\Lambda_0^i)\geq R_B$, set $j=j+1$ and go to step 3; otherwise, set $(r,\lambda,\theta_P)=\Pi(\mathcal{T}^\mu(t,\Theta_0^i,\Lambda_0^i))$ and $u=U^\mu(\bar t,\Theta_0^i,\Lambda_0^i)$, and go to step 6.	
			\item[6.] Set $\mathcal{D} = \mathcal{D}\cup \{[r,\lambda,\theta_P,\mu^j,u]\}$, set $\bar t = \bar t + h $, and go to step 4.	
			\item[7.] End.
		\end{itemize}
		%\DrawLine
	\end{tcolorbox}
\end{center}

By {\it Procedure \ref{pro1}}, the dataset for the mapping $(r,\lambda,\theta_P,\mu)\mapsto C(r,\lambda,\theta_P,\mu)$ is eventually included in the set $\mathcal{D}$. An NN can be trained by the dataset $\mathcal{D}$ to approximate the mapping $(r,\lambda,\theta_{P},\mu)\mapsto C(r,\lambda,\theta_{P},\mu)$. Set the numbers $N$, $M$, $h$, and $R_B$ in {\it Procedure \ref{pro1}} as $1000$, $10$, $0.5$ sec, and $30$ km, respectively. This means that $1\times10^4$ trajectories are generated by {\it Procedure \ref{pro1}}. Then, a Feedforward NN (FNN) comprising three hidden layers, each with 20 neurons, is trained by the dataset to approximate the optimal feedback control $C(r,\lambda,\theta_P,\mu)$. The loss function is chosen as the mean-squared error between the predicted outputs and the actual values within dataset $\mathcal{D}$. Finally, the training is terminated when the mean-squared error reaches $1\times10^{-5}$. Let $N(r,\lambda,\theta_P,\mu)$ be the trained FNN. It is capable of computing an optimal guidance command in approximately 0.16 milliseconds for any valid input  $(r,\lambda,\theta_P,\mu)$. This inference speed is achieved on a platform with MYC-Y6ULY2 CPU at $528$ MHz. 

It should be noted that the trained network $N(r,\lambda,\theta_P,\mu)$ cannot be directly used once $r>R_B$ or the heading angle of the target is not zero. According to Lemma \ref{lemma:rotation}, for any  $(r,\lambda,\theta_P,\theta_T,\mu)$ with $\theta_T \neq 0$, we can use $N(r,\lambda+\theta_T,\theta_P+\theta_T,\mu)$ to approximate the optimal feedback control  $u^*(r,\lambda,\theta_P,\theta_T,\mu)$. According to Lemma \ref{lemma5}, if $r>R_B$, we can use $N(r/s,\lambda+\theta_T,\theta_P+\theta_T,\mu)/s$ to approximate the optimal feedback control  $u^*(r,\lambda,\theta_P,\theta_T,\mu)$. In the subsequent section, the trained neural network will also be modified for intercepting maneuvering targets.

\section{Guidance for Intercepting Maneuvering Targets}\label{sec:GMT}
%In this section, we first present a guidance law with a bias term devised from the perspective of relative motion to compensate the effect of target's maneuver. By the compensation of the bias term, we just need to generate the dataset for intercepting nonmaneuvering target. Then, we can train a simple Feedforward Neural Network (FNN) to approximate the relationship between the relative state and optimal control command.

\begin{figure*}[!htp]
	\centering
	\begin{subfigure}{0.4\linewidth}
		\centering
		\includegraphics[width = 0.99\linewidth]{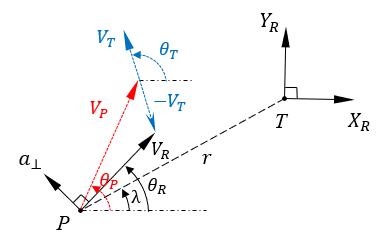}
		\caption{The relative frame}
		\label{subfig:vr}
	\end{subfigure}
	~~~~~
	\begin{subfigure}{0.4\linewidth}
		\centering
		\includegraphics[width = 0.99\linewidth]{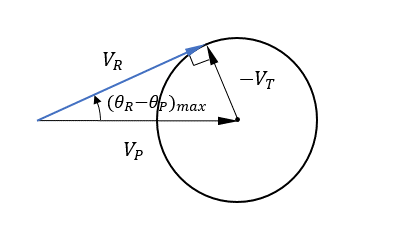}
		\caption{The relations between $V_R$ and $V_P$}
		\label{subfig:vr2}
	\end{subfigure}
	\caption{The geometry for relative motion.}
	\label{fig:vr}
\end{figure*}

Note that the conventional PN is often modified to design guidance laws for intercepting maneuvering targets in the literature; see, e.g., \cite{zarchan2012tactical,7472963,jeon2015exact}. In general, the PN-like guidance command for intercepting maneuvering targets takes a form of 
\begin{align}\label{eq:a_p}
	u = F(t) u_{PN} + G(t) a_T/V_P
\end{align}
where $u_{PN}$ is the PN-like guidance command, $a_T$ is the lateral acceleration of the target, $F(t)$ is the time varying gain of the PN term, and $G(t)$ is the time varying gain of the biased term. Three typical guidance laws, taking the form in Eq.~(\ref{eq:a_p}), for intercepting maneuvering targets are listed in Table \ref{tab:gui_form}, where $V_R$ is the relative speed, and $\theta_R$ is the relative heading angle, as shown in Fig.~\ref{fig:vr}.
\begin{table}[hbt!]
	\caption{\label{tab:gui_form} Three typical guidance laws taking the form of Eq. (\ref{eq:a_p}).}
	\centering
	\renewcommand{\arraystretch}{2}
	\begin{tabularx}{0.55\textwidth}{XXXX}
		\toprule
		& $u_{PN}$    & $F(t)$    & $G(t)$\\ \midrule
		APN \cite{zarchan2012tactical} & $-3\dot r\dot{\lambda}$ & $\displaystyle{\frac{1}{\cos\sigma_P}}$ &  $\displaystyle{\frac{2\sin\sigma_T}{3\cos\sigma_P}}$\\
		SMG \cite{7472963}  & $(\frac{800}{V_P}-2\dot r)\dot{\lambda}$ & $\displaystyle{\frac{1}{\cos\sigma_P}}$ & $\displaystyle{\frac{\sin\sigma_T}{\cos\sigma_P}}$\\
		PEOG \cite{jeon2015exact} & $3V_R\dot{\lambda}/V_P$ & $\displaystyle{\frac{1}{\cos(\theta_P-\theta_R)}}$& $\displaystyle{\frac{\cos(\theta_T-\theta_R)}{\cos(\theta_P-\theta_R)}}$\\ \bottomrule
	\end{tabularx}
\end{table}

It is apparent from Table \ref{tab:gui_form} that the APN and the SMG are both singular if $\sigma_P = \pi/2$. To address this singularity, the PEOG was designed by using the optimal control theory. According to Fig.~\ref{subfig:vr2}, we have that $\cos(\theta_R - \theta_P)$ cannot be zero if $\mu < 1$, indicating that the PEOG is not singular even if $\sigma_P = \pi/2$. Since we have obtained the nonlinear optimal guidance command in the previous sections, replacing the term $u_{PN}$ of the PEOG in Table \ref{tab:gui_form} with our nonlinear optimal guidance command should be able to further improve the performance for intercepting maneuvering targets. By replacing the term $3 V_R \dot\lambda/V_P$ of the PEOG in Table \ref{tab:gui_form} with the nonlinear optimal guidance command $N(r,\lambda+\theta_T,\theta_P+\theta_T,\mu)$, we accordingly propose the following guidance law for intercepting maneuvering targets:
\begin{equation}\label{eq:mod_peog}
	u=\displaystyle{\frac{ N(r,\lambda+\theta_T,\theta_P+\theta_T,\mu)}{\cos(\theta_P-\theta_R)}}+\displaystyle{\frac{\cos(\theta_T-\theta_R)a_T}{\cos(\theta_P-\theta_R)V_P}}.
\end{equation}

In the following section, we shall show by numerical examples that the guidance law in Eq.~(\ref{eq:mod_peog}) performs better than the existing PEOG in Table \ref{tab:gui_form}.

\section{Numerical Simulations}\label{Sec:NS}
%Thus, the trained FNN is able to embedded into a guidance system which frequency is more than $3$ KHz. Note that the computational period is more than enough for the onboard computational platform of pursuer.
Three scenarios will be simulated to demonstrate the developments of the paper. 

\subsection{Simulations for Intercepting Nonmaneuvering Target}\label{subsecA}
In this subsection, two engagements are simulated to compare the developed nonlinear optimal guidance with existing guidance laws, and for simplicity we shall use NOG to denote the developed nonlinear optimal guidance. 
\subsubsection{Engagement I} %Different Initial Heading Angles for the Pursuer with a Fixed Heading Angle for the Target
Consider an initial position of the pursuer as  $(3000,4000)$ m. Five different initial heading angles, $165~\deg$, $195~\deg$, $225~\deg$, $255~\deg$, and $285~\deg$, for the pursuer are considered for illustration.  Initially, the target is located at the origin $(0,0)$ m with its heading angle fixed as $180$ deg. The speed of the pursuer is set as $V_P = 1000$ $\mathrm{m/s}$, and the speed of the target is set as $V_T = 400$ $\mathrm{m/s}$. The acceleration of the pursuer is considered to be limited within $50$ g, where g $ = 9.81$ $\mathrm{m/s^2}$ is the gravitational acceleration constant at sea level. Then, the trained FNN is used to generate the NOG command, as shown by the closed-loop diagram in Fig.~\ref{fig:Guidance}. The trajectories related to NOG are represented by the solid lines in Fig.~\ref{fig:traj-case1}.

The blue dotted-dashed lines in Fig.~\ref{fig:traj-case1} denote the trajectories related to the PEOG and the trajectories generated by PN are presented by the black dashed lines in Fig.~\ref{fig:traj-case1}. The corresponding guidance command profiles of different guidance laws are reported in Fig.~\ref{fig:a-case1}, and the profiles of pursuer's heading angle are shown in Fig.~\ref{fig:th-case1}. Denote by $\bar J=1/2\int_{0}^{t_f}a_P^2\ dt$ the control effort required by the pursuer for intercepting the target. Define $\theta_I$ as the ideally initial heading angle of the pursuer that can achieve the zero-effort collision triangle \cite{cho2014optimal}. To compare the optimality of these guidance laws, the values of $\bar J$ are presented in Table~\ref{tab:table1}.

\begin{table*}[hbt!]
	\caption{\label{tab:table1} Engagement I: The values of control effort required by the pursuer for intercepting the target.}
	\centering
	\begin{tabular}{ccccc}
		\toprule
		$|\theta_{P0}-\theta_I|$&$\theta_{P0}$&\multicolumn{3}{c}{$\bar J$ ($\mathrm{m^2/s^3}$)}\\\cline{3-5}
		(deg) & (deg) & NOG & PN & PEOG\\\midrule
		\multicolumn{1}{c}{10.53} & 225 & $7.2408\times10^3$& $7.8048\times10^3$& $7.2684\times10^3$\\%-135
		\multicolumn{1}{c}{19.47} & 195 & $2.3522\times10^4$& $2.4819\times10^4$& $2.3871\times10^4$\\%-165
		\multicolumn{1}{c}{40.53} & 255 & $1.0269\times10^5$& $1.1516\times10^5$& $1.1103\times10^5$\\%-105
		\multicolumn{1}{c}{49.47} & 165 & $1.3085\times10^5$& $1.3710\times10^5$& $1.3910\times10^5$\\%-195
		\multicolumn{1}{c}{70.53} & 285 & $2.6008\times10^5$& $3.0527\times10^5$& $2.9458\times10^5$\\%-75
		\bottomrule
	\end{tabular}
\end{table*}

As we can see from Fig.~\ref{fig:traj-case1}, the trajectories related to different guidance laws are coincident when the initial heading angle of pursuer is close to the collision course. In the meantime, the values of control effort are also close to each other, as shown in Table~\ref{tab:table1}. Nevertheless, it is clear to distinguish the trajectories of different guidance laws if the initial heading angle of the pursuer is far from the collision course. We can also see from Table~\ref{tab:table1} that the trajectories by NOG is better than PN and PEOG in terms of control efforts. In fact, it can be seen that the farther the initial pursuer's heading angle is from the collision course, the more different of the control effort is.

\subsubsection{Engagement II}%: Different Heading Angles for the Target with A Fixed Initial Heading Angle for the Pursuer}
\begin{figure*}[!htp]
	\centering
	\begin{subfigure}{0.46\linewidth}
		\centering	
		\includegraphics[width = 8cm]{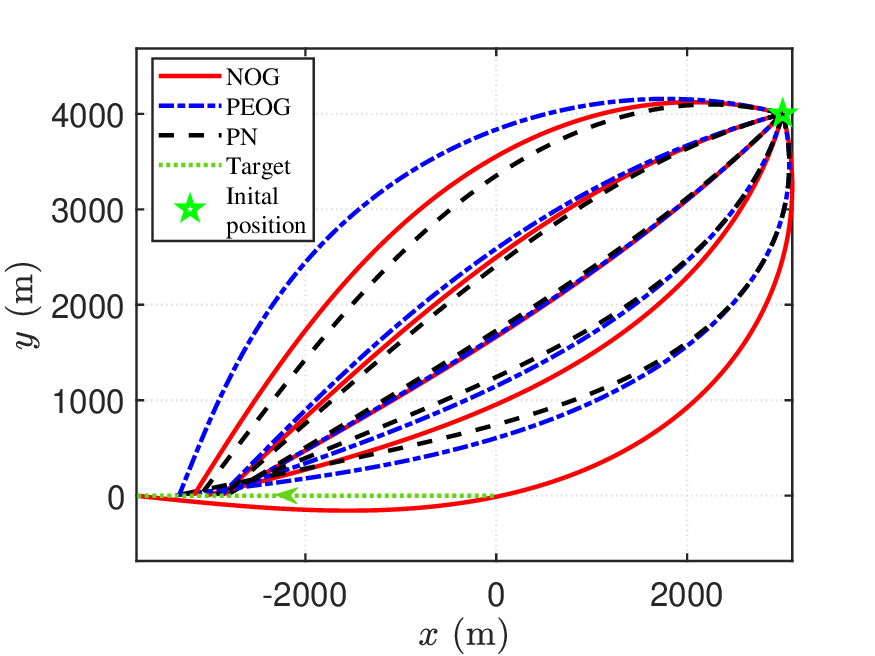}
		\caption{Engagement I: Trajectories related to different guidance laws.}
		\label{fig:traj-case1}
	\end{subfigure}
	~~~~~
	\begin{subfigure}{0.46\linewidth}
		\centering	
		\includegraphics[width = 8cm]{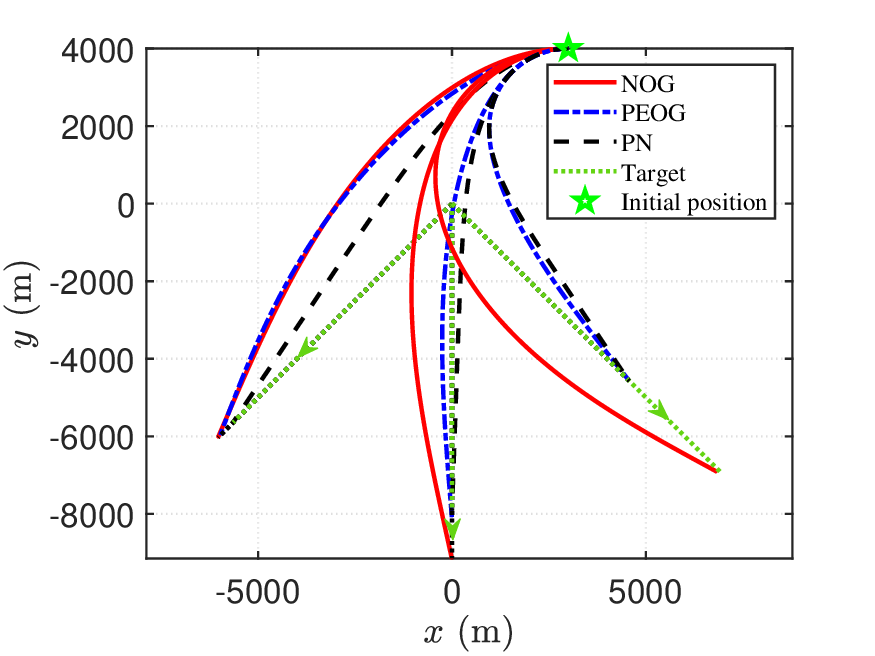}
		\caption{Engagement II: Trajectories related to different guidance laws.}
		\label{fig:traj-case2}
	\end{subfigure}

	\begin{subfigure}{0.46\linewidth}
		\centering	
		\includegraphics[width = 8cm]{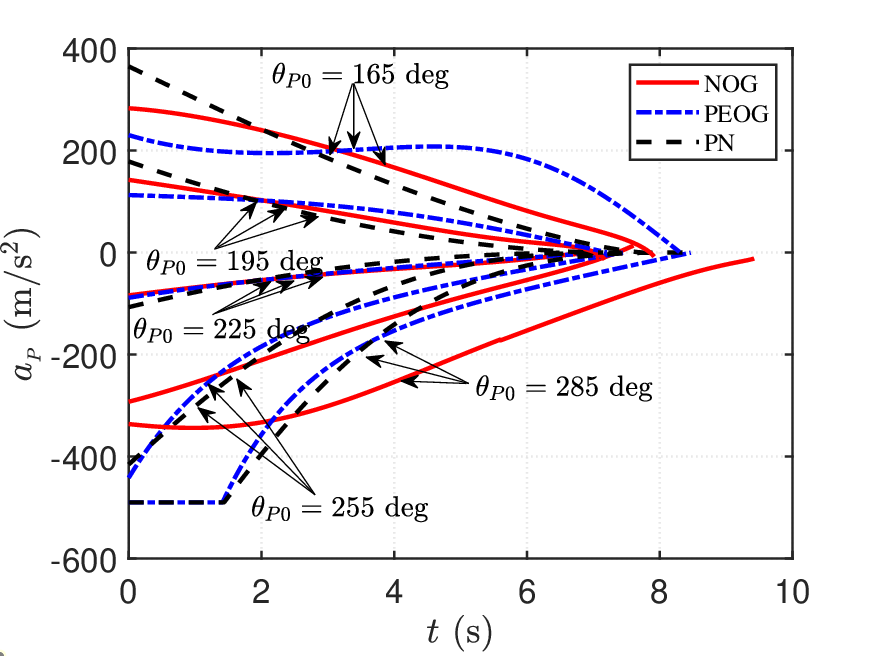}
		\caption{Engagement I: Guidance command profiles of different guidance laws.}
		\label{fig:a-case1}
	\end{subfigure}
	~~~~~
	\begin{subfigure}{0.46\linewidth}
		\centering	
		\includegraphics[width = 8cm]{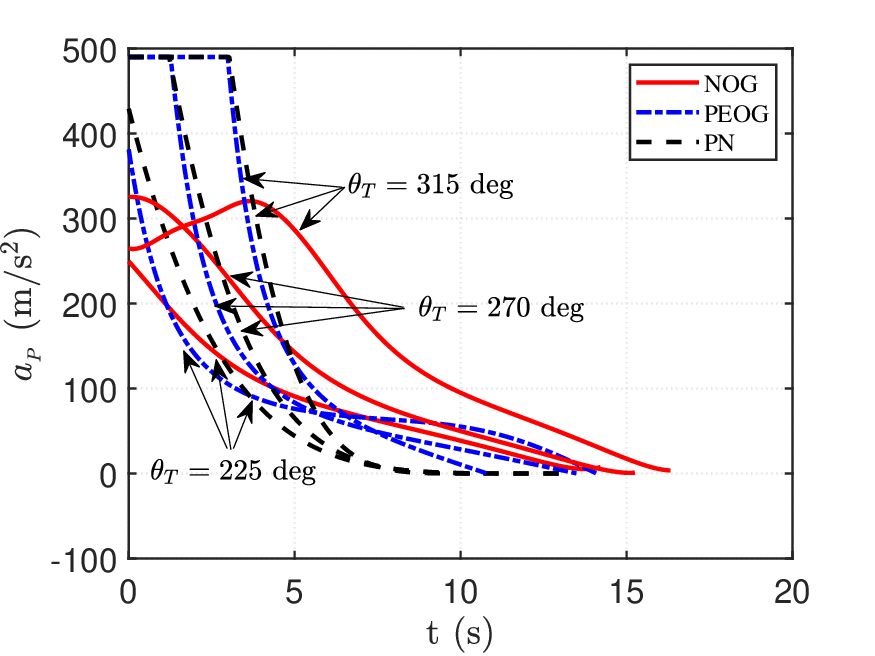}
		\caption{Engagement II: Guidance command profiles of different guidance laws.}
		\label{fig:a-case2}
	\end{subfigure}

	\begin{subfigure}{0.46\linewidth}
		\centering	
		\includegraphics[width = 8cm]{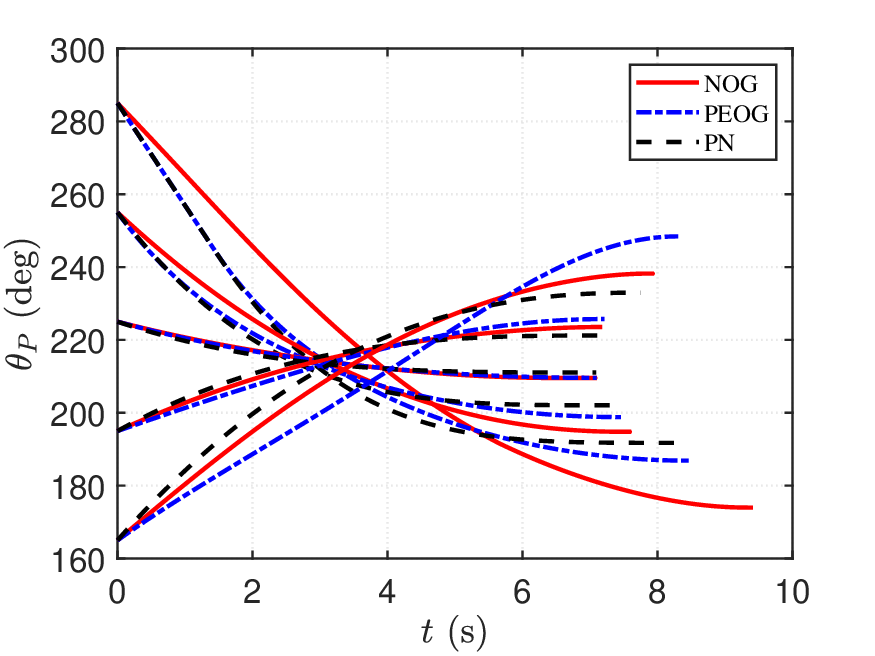}
		\caption{Engagement I: Profiles of pursuer's heading angle related to different guidance laws.}
		\label{fig:th-case1}
	\end{subfigure}
	~~~~~
	\begin{subfigure}{0.46\linewidth}
		\centering	
		\includegraphics[width = 8cm]{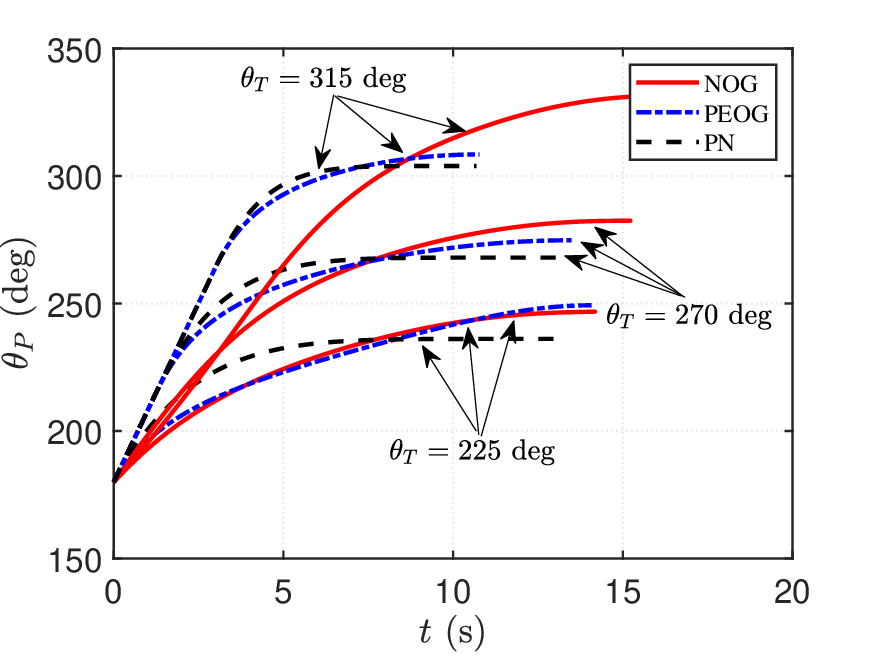}
		\caption{Engagement II: Profiles of pursuer's heading angle related to different guidance laws.}
		\label{fig:th-case2}
	\end{subfigure}
		\caption{Trajectories, profiles of guidance command and pursuer's heading angle of Engagement I \& II.}
		\label{fig:case2}
\end{figure*}
For Engagement II, the initial state of the pursuer is considered to be
\begin{equation*}
	(x_{P0},y_{P0},\theta_{P0})=(3000\ \rm{m},4000\ \rm{m},180\ \rm{deg}).
\end{equation*}
Set the speed of the target as $600~\mathrm{m/s}$. The target is initially located at the origin $(0,0)$ m. Three different heading angles for the target are considered for illustration. The other simulation parameters are the same as those in Engagement I. The trajectories related to the NOG, the PN, and the PEOG with different target's heading angles are presented in Fig.~\ref{fig:traj-case2}. The corresponding guidance command profiles of three guidance laws are reported in Fig.~\ref{fig:a-case2}. The profiles of pursuer's heading angles are depicted in Fig.~\ref{fig:th-case2}. In addition, the values of control effort consumed for different target's heading angles are reported in Table~\ref{tab:table2}. Similar to the result in Engagement I, the values of control effort required by the NOG are much lower than those related to the PN and the PEOG, as shown in Table \ref{tab:table2}. As depicted in Fig.~\ref{fig:a-case2}, the NOG has a lower requirement on the normal acceleration.

\begin{table*}[hbt!]
	\caption{\label{tab:table2} Engagement II: The values of control effort required by the pursuer for intercepting the target.}
	\centering
	\begin{tabular}{ccccc}
		\toprule
		$|\theta_{P0}-\theta_I|$&$\theta_{T0}$&\multicolumn{3}{c}{$\bar J$, $\mathrm{m^2/s^3}$}\\\cline{3-5}
		(deg) & (deg) & NOG & PN & PEOG\\\midrule
		\multicolumn{1}{c}{49.89}& 225& $7.6826\times10^4$& $1.1520\times10^5$ & $8.5244\times10^4$\\%-135
		\multicolumn{1}{c}{69.87}& 270 & $1.8438\times10^5$& $2.7319\times10^5$ & $2.5202\times10^5$\\%-90
		\multicolumn{1}{c}{76.46}& 315 & $3.1340\times10^5$ &$4.5630\times10^5$ & $4.4117\times10^5$\\%-45
		\bottomrule
	\end{tabular}
\end{table*}

\subsection{Simulations for Intercepting Constant Maneuvering Target}\label{subsecB}

In this subsection, we consider to intercept a constant maneuvering target. The lateral acceleration of the target is set as $-50\ \mathrm{m/s^2}$. The initial location and the initial heading angle of the target are $(0,0)$ m and $180$ deg, respectively. The pursuer is initially located at $(10000,-10000)$ m, and the initial heading angle is $180$ deg. Set the speed of the target as $500~\mathrm{m/s}$. The other simulation parameters are the same as those in Engagement I of Subsection \ref{subsecA}. The improved NOG in Eq.~(\ref{eq:mod_peog}) and the NOG are all used in this example to compare with existing guidance laws. 

\begin{figure}[hbt!]
	\centering	
	\includegraphics[width=.5\textwidth]{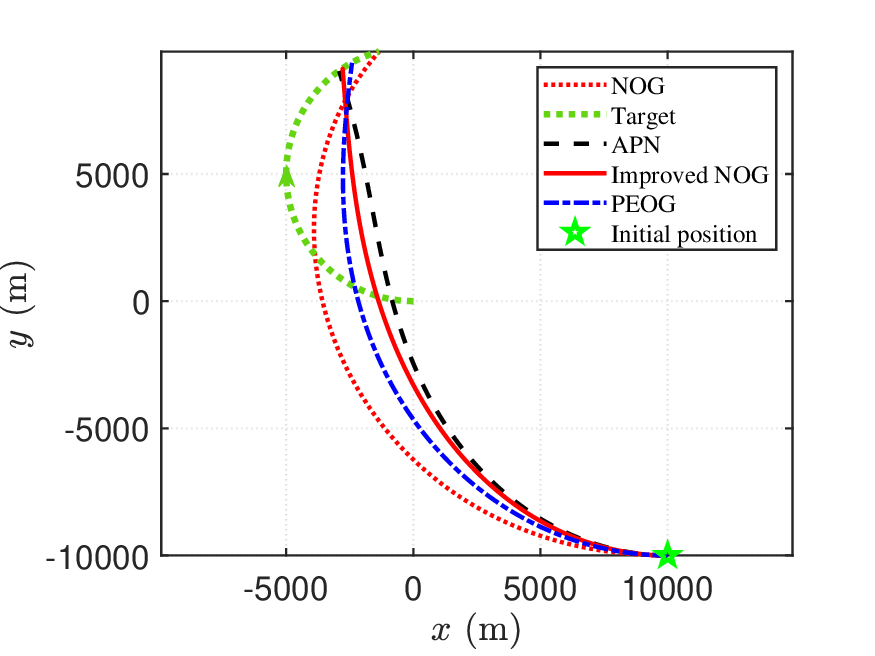}
	\caption{Comparison of trajectories related to different guidance laws for intercepting a constant maneuvering target.}
	\label{fig:traj-case3}
\end{figure}

\begin{figure}[hbt!]
\centering	
\includegraphics[width=.5\textwidth]{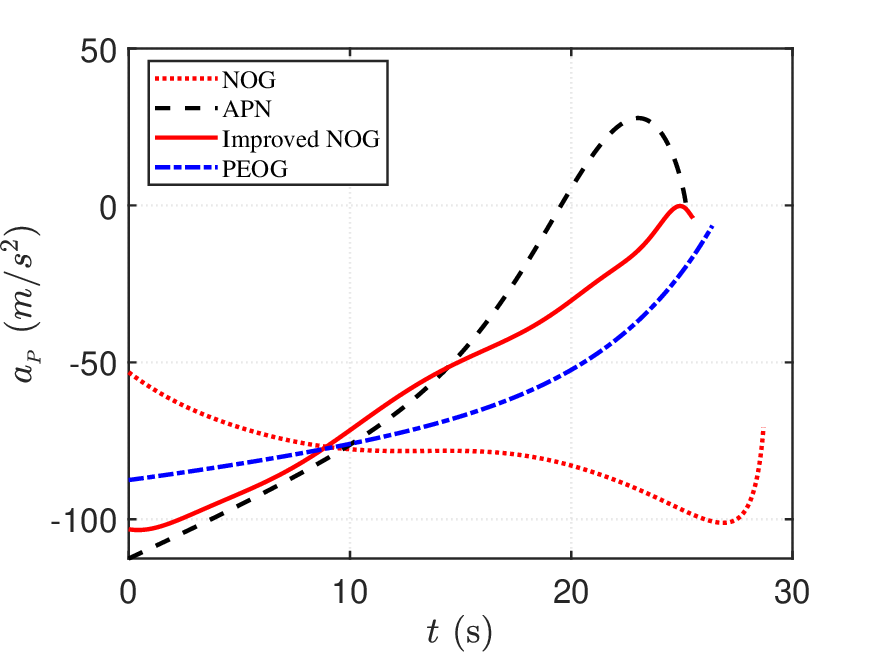}
\caption{Guidance command profiles of different guidance laws for intercepting a constant maneuvering target.}
\label{fig:a-case3}
\end{figure}

\begin{figure}[hbt!]
\centering	
\includegraphics[width=.5\textwidth]{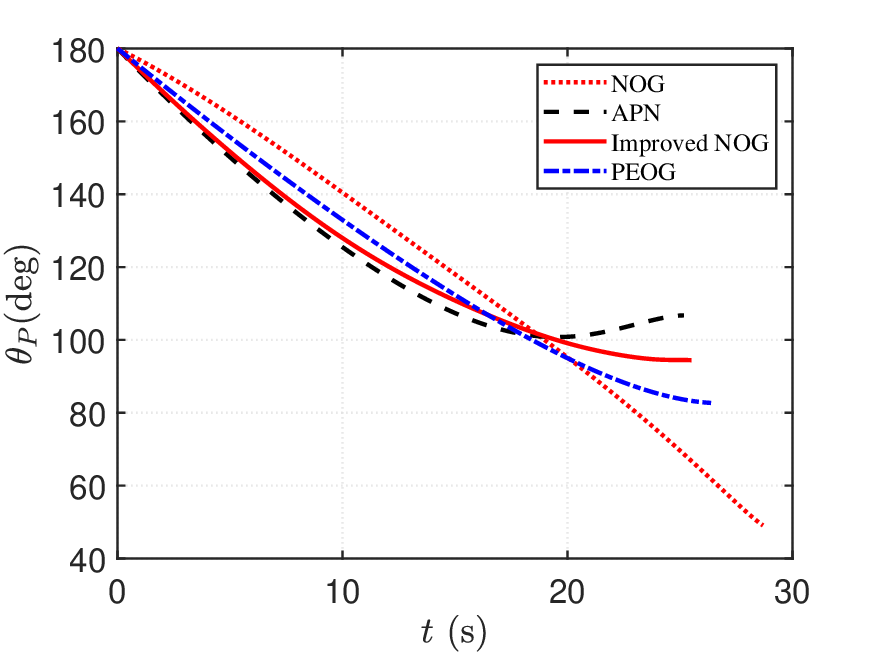}
\caption{Profiles of pursuer's heading angle related to different guidance laws.}
\label{fig:th-case3}
\end{figure}

%\begin{table}[hbt!]
%	\caption{\label{tab:table3} The values of control effort consumed by different guidance laws.}
%	\centering
%	\begin{tabular}{ccccccc}
%		\toprule
%		  &$\bar J$ of APN, $m^2/s^3$& &$\bar J$ of PEOG, $m^2/s^3$& &$\bar J$ of SMG, $m^2/s^3$\\\midrule
%		\multicolumn{1}{c}{Origin Method}& $2.2744\times10^4$ && $1.6966\times10^5$&& $2.1126\times10^5$\\
%		\multicolumn{1}{c}{FNN Modified}& $1.9119\times10^4$ && $1.5493\times10^5$ && $1.7047\times10^5$\\
%		\bottomrule
%	\end{tabular}
%\end{table}
The trajectory related to the NOG is shown by the dotted curve in Fig.~\ref{fig:traj-case3}.  The trajectory generated by the improved NOG is represented by the solid curve in Fig.~\ref{fig:traj-case3}. The trajectories generated by the APN and PEOG  are respectively indicated by the dashed curve and the dotted-dashed curve in Fig.~\ref{fig:traj-case3}. The control effort related to the APN is $5.9004\times10^4\ \mathrm{m^2/s^3}$, and that related to PEOG is $6.0079\times10^4\ \mathrm{m^2/s^3}$. Without the augmented information of the target's maneuver, the control effort required by the NOG is $9.2690\times10^4\ \mathrm{m^2/s^3}$. However, the control effort required by the  improved NOG is reduced to $5.5851\times10^4\ \mathrm{m^2/s^3}$. The control profiles and the time histories for pursuer's heading angles are presented in Fig.~\ref{fig:a-case3} and Fig.~\ref{fig:th-case3}, respectively. Notice from Fig.~\ref{fig:a-case3} that the absolute value of normal acceleration required by the improved NOG at the terminal phase is smaller than that by the NOG.

\subsection{Simulations for Intercepting Variable Maneuvering Target}\label{subsecC}
In this scenario, the target's lateral acceleration is chosen as $100\sin0.5t$ $\mathrm{m/s^2}$. The initial states of the pursuer and the target are set as
\begin{equation*}
	\begin{split}
		(x_{P0},y_{P0},\theta_{P0})&=(-10000\ \mathrm{m},-8000\ \mathrm{m},150\ \deg),\\
		(x_{T0},y_{T0},\theta_{T0})&=(0\ \mathrm{m},0\ \mathrm{m},120\ \deg).
	\end{split}	
\end{equation*}
Set the speed of the target as $700~\mathrm{m/s}$. The other simulation parameters are the same as Engagement I in Subsection \ref{subsecA}. In order to illustrate the guidance performance of intercepting variable maneuvering target, we consider the autopilot dynamics as a first order lag system, whose time constant is chosen as $0.5$ s \cite{doi:10.2514/6.2009-6089}.

\begin{figure}[hbt!]
	\centering	
	\includegraphics[width=.5\textwidth]{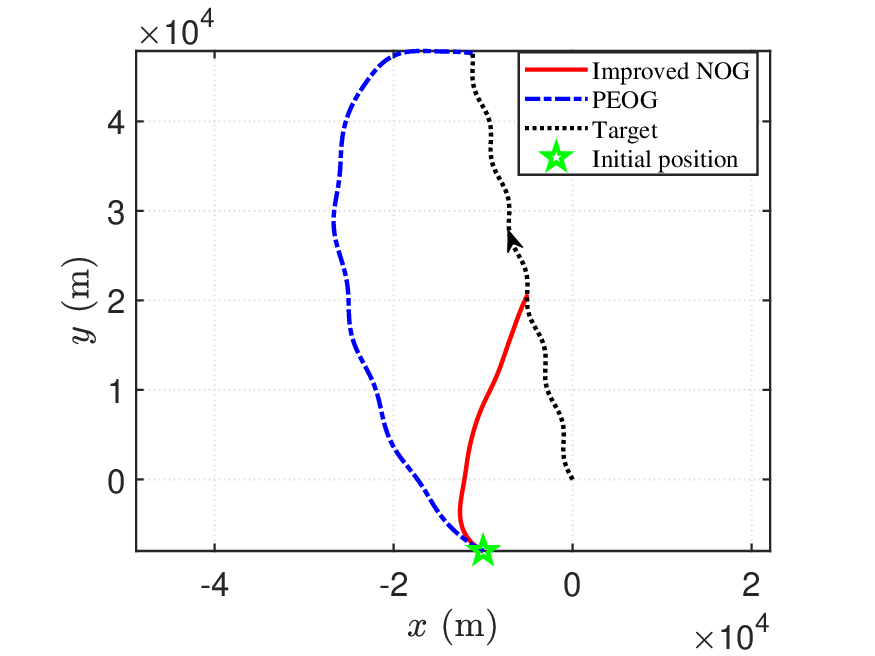}
	\caption{Comparison of trajectories related to different guidance laws for intercepting a variable maneuvering target.}
	\label{fig:traj-case4}
\end{figure}

\begin{figure}[hbt!]
	\centering	
	\includegraphics[width=.5\textwidth]{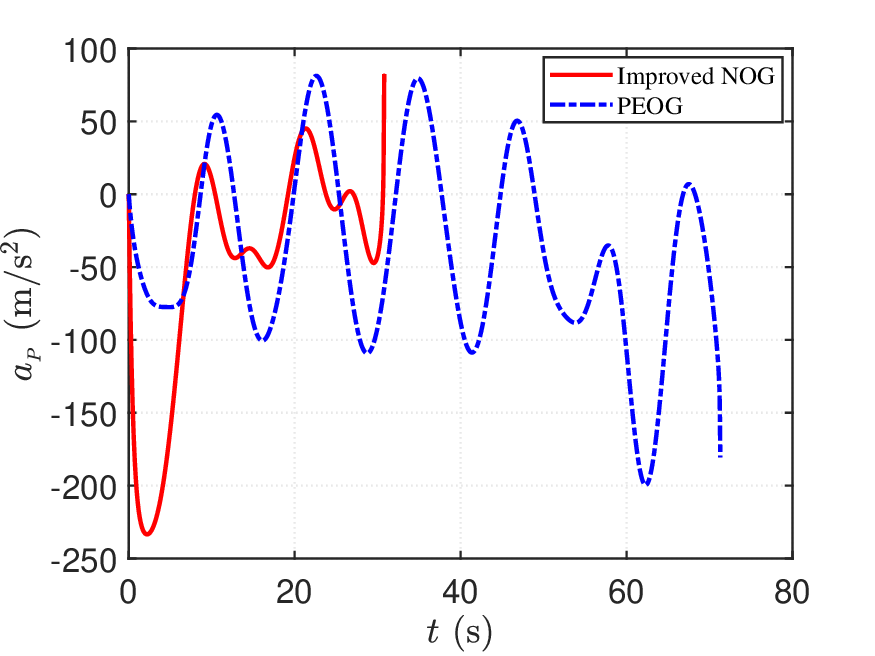}
	\caption{Guidance command profiles of different guidance laws for intercepting variable a maneuvering target.}
	\label{fig:a-case4}
\end{figure}

\begin{figure}[hbt!]
	\centering	
	\includegraphics[width=.5\textwidth]{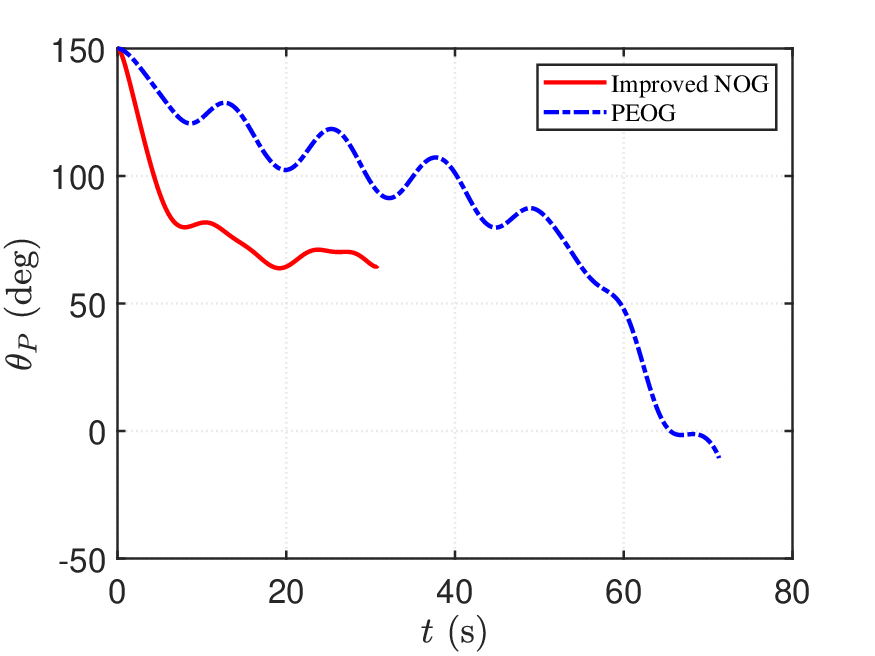}
	\caption{Profiles of pursuer's heading angle related to different guidance laws.}
	\label{fig:th-case4}
\end{figure}

%Denote by $a_{Pf}$ the terminal value of the pursuer's normal acceleration. $a_{Pf}$ and $\bar J$ of those guidance laws are compared in the Table.~\ref{tab:table4}. 
%\begin{table}[hbt!]
%	\caption{\label{tab:table4} The values of $a_{Pf}$ and $\bar J$ consumed by different guidance laws.}
%	\centering
%	\begin{tabular}{ccccccccc}
%		\toprule
%		&\multicolumn{2}{c}{APN}& & \multicolumn{2}{c}{PEOG}& & \multicolumn{2}{c}{SMG}\\
%		\cline{2-3} \cline{5-6} \cline{8-9}
%		& $a_{Pf}$, $m/s^2$ & $\bar J$, $m^2/s^3$ &&$a_{Pf}$, $m/s^2$ & $\bar J$, $m^2/s^3$ && $a_{Pf}$, $m/s^2$ & $\bar J$, $m^2/s^3$\\\midrule
%		\multicolumn{1}{c}{Origin Method}& $19.2$& $2.0435\times10^5$ && $29.3$& $2.0091\times10^5$&& $42.3$& $2.0412\times10^5$\\
%		\multicolumn{1}{c}{FNN Modified}& $1.02$& $1.4493\times10^5$ && $-1.36$& $1.2802\times10^5$ && $20.8$& $1.6716\times10^5$\\
%		\bottomrule
%	\end{tabular}
%\end{table}
Because the initial lead angle of the pursuer is close to $\pi/2$, the APN guidance commands diverged at the beginning. Thus, the results are compared with the conventional PN and the PEOG. Due to the presence of autopilot lag, the pursuers guided by the PN and the NOG miss the target. Thus, we present the trajectories of the PEOG and the improved NOG only. The trajectory of the improved NOG is shown by the solid curve in Fig.~\ref{fig:traj-case4}. The trajectory of the PEOG is shown by the the dotted-dashed curve in Fig.~\ref{fig:traj-case4}.
 The total control effort required by the improved NOG is $1.2683\times10^5\ \mathrm{m^2/s^3}$, which is significantly lower than $1.9897\times10^5\ \mathrm{m^2/s^3}$ required by the PEOG. The profiles of guidance command and pursuer's heading angle are presented in Fig.~\ref{fig:a-case4} and Fig.~\ref{fig:th-case4}, respectively.

\section{Conclusions}\label{Sec:Con}
This paper aims to address the nonlinear optimal guidance problem for intercepting moving targets. The core objective is to develop a real-time optimal feedback control strategy for guiding a pursuer to intercept a moving target. The necessary optimality conditions for the corresponding nonlinear optimal control problem were derived by using PMP, and these necessary conditions were further employed to show that the extremal trajectories and extremal controls are determined by two scalars. In addition, two extra optimality conditions were established to ensure local optimality. By analyzing the geometric properties of extremal trajectories, some properties of the optimal feedback control were studied, and were embedded into the procedure for generating the dataset for the mapping from state to optimal feedback control. Thus, the size of the dataset was significantly reduced. This allows to train a lightweight FNN to approximate the optimal feedback control in real time. Moreover, according to the existing guidance laws for intercepting maneuvering targets, the trained FNN was adjusted so that an FNN-based guidance law for intercepting maneuvering targets was proposed. Numerical simulations show that the proposed nonlinear optimal guidance outperforms the existing guidance laws.  

\section*{Appendix}

\subsection{Proof of Lemma \ref{lemma:ex_tra}}\label{app:A}
\setcounter{equation}{0}
\renewcommand{\theequation}{A\arabic{equation}}
 According to PMP, for any optimal trajectory $(r(t),\lambda(t),\theta_P(t))$, $t\in [0,t_f]$, of the OCP with $\mu \in (0,1)$, there exists a continuous costate function $(p_r(t),p_\lambda(t),p_\theta(t))\in\mathbb{R}^3$ for $t\in [0,t_f]$ so that $(r(t),\lambda(t),\theta_P(t),p_r(t),p_\lambda(t),p_\theta(t))$ is a solution of the canonical equations combining Eq.~(\ref{eq:norm-r-lambda}) and Eq.~(\ref{eq:ode1}). By the definition of differential equations in Eq.~(\ref{equation:phpx1}), for any $\mu\in(0,1)$, a trajectory $(R^\mu(\tau,\Theta_0,\Lambda_0),\Lambda^\mu(\tau,\Theta_0,\Lambda_0),\Theta^\mu(\tau,\Theta_0,\Lambda_0))$ for $\tau\in[0,t_f]$ can be obtained by solving the initial value problem defined in Eq.~(\ref{equation:phpx1}) and Eq.~(\ref{eq:H=0}) with $\Theta_0\in[0,2\pi)$ and $\Lambda_0\in[0,2\pi)$ so that
\begin{equation*}
	\begin{split}
		(R^\mu(\tau,\Theta_0,\Lambda_0),\Lambda^\mu(\tau,\Theta_0,\Lambda_0),\Theta^\mu(\tau,\Theta_0,\Lambda_0))=(r(t_f-\tau),\lambda(t_f-\tau),\theta_P(t_f-\tau))
	\end{split}
\end{equation*}
completing the proof of the first statement of Lemma \ref{lemma:ex_tra}.

Given any $\mu \in (0,1)$, and any $\Theta_0$ and $\Lambda_0$, set $\tau = t_f-t$, and let
\begin{equation*}%\label{eq:A1}
	\begin{cases}
		r(t)=R^\mu(\tau,\Theta_0,\Lambda_0)\\
		\lambda(t) = \Lambda^\mu(\tau,\Theta_0,\Lambda_0)\\
		\theta_P(t)= \Theta^\mu(\tau,\Theta_0,\Lambda_0)\\		
	\end{cases}
\end{equation*}
By the initial value problem defined in Eq.~(\ref{equation:phpx1}) and Eq.~(\ref{eq:H=0}), the trajectory $(r(t),\lambda(t),\theta_P(t))$ for $t\in [0,t_f]$ satisfies the necessary conditions in Eqs.~(\ref{sample:Ham}-\ref{equation:Hamiltonian}). Thus, the trajectory $(r(t),\lambda(t),\theta_P(t))$ for $t\in [0,t_f]$ is an extremal trajectory, completing the proof of the second statement of Lemma \ref{lemma:ex_tra}.

%By contradiction, assume that for any extremal trajectory $(r(t),\lambda(t),\theta_P(t))$, $t\in [0,t_f]$, there exists more than one pairs $(\Theta_0,\Lambda_0)$ so that
%\begin{equation}\label{eq:pro_lm1}
%	(R^\mu(\tau,\Theta_0,\Lambda_0),\Lambda^\mu(\tau,\Theta_0,\Lambda_0),\Theta^\mu(\tau,\Theta_0,\Lambda_0))=(r(t_f-\tau),\lambda(t_f-\tau),\theta_P(t_f-\tau))
%\end{equation}
%According to the proof of the existence, different trajectories $(R^\mu(\tau,\Theta_0,\Lambda_0),\Lambda^\mu(\tau,\Theta_0,\Lambda_0),\Theta^\mu(\tau,\Theta_0,\Lambda_0))$ for $\tau\in[0,t_f]$ can be obtained by solving the initial value problem defined in (\ref{equation:phpx1}) and (\ref{eq:H=0}) with pairs of $(\Theta_0,\Lambda_0)$, which contradicts with the condition of (\ref{eq:pro_lm1}). Thus, there is only one pair $(\Theta_0,\Lambda_0)$ so that $(R^\mu(\tau,\Theta_0,\Lambda_0),\Lambda^\mu(\tau,\Theta_0,\Lambda_0),\Theta^\mu(\tau,\Theta_0,\Lambda_0))=(r(t_f-\tau),\lambda(t_f-\tau),\theta_P(t_f-\tau))$, completing the proof.
$\hfill\square$

\subsection{Proof of Lemma \ref{lemma:3}}\label{app:B}
\setcounter{equation}{0}
\setcounter{figure}{0}
\renewcommand{\theequation}{B\arabic{equation}}
\renewcommand{\thefigure}{B\arabic{figure}}
By contradiction, assume that along the optimal trajectory $\Pi(\mathcal{T}^\mu(\cdot,\Theta_0,\Lambda_0))$ on $[0,\tau_f]$, there exists a time $\bar\tau$ within $(0,\tau_f)$ such that the vector $[\cos\Theta^\mu(\bar\tau,\Theta_0,\Lambda_0),\sin\Theta^\mu(\bar\tau,\Theta_0,\Lambda_0)]$ is collinear with the LOS; i.e., Eq.~(\ref{eq:sin_dth}) holds. Let A be the state at $\bar\tau$; i.e., $A=\Pi(\mathcal{T}^\mu(\bar\tau,\Theta_0,\Lambda_0))$. Without loss of generality, assume that the initial location and the heading angle of the target are $(0,0)$ and $0$ deg, respectively. Then, we can obtain the extremal trajectory $\Pi(\mathcal{T}^\mu(\tau,\Theta_0,\Lambda_0))$ in the OXY frame by the following equations
\begin{equation*}
	\begin{split}
		x(\tau) &= R^\mu(\tau,\Theta_0,\Lambda_0)\cos\Lambda^\mu(\tau,\Theta_0,\Lambda_0)+\mu\tau\\
		y(\tau) &= R^\mu(\tau,\Theta_0,\Lambda_0)\sin\Lambda^\mu(\tau,\Theta_0,\Lambda_0)\\
	\end{split}
\end{equation*}

\begin{figure}[hbt!]
	\centering	
	\includegraphics[width=.4\textwidth]{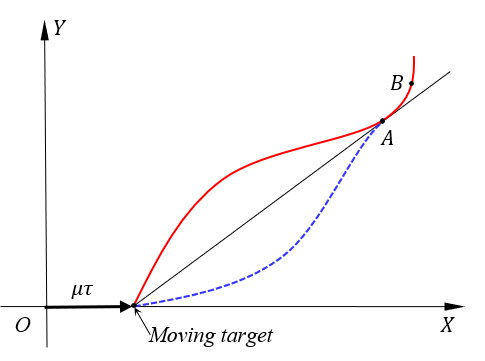}
	\caption{Geometry for existing points with $\Theta=\Lambda$.}
	\label{fig:A1}
\end{figure}

Set $\hat U(\tau,\Theta_0,\Lambda_0)=-U^\mu(\tau,\Theta_0,\Lambda_0)$. Then, we have that $\hat U(\tau,\Theta_0,\Lambda_0)$ for $\tau\in[0,\bar\tau]$ is the control for the pursuer to move along the symmetric path from A to the target, as shown by the dashed curve in Fig.~\ref{fig:A1}, where the solid curve denotes the extremal trajectory $\Pi(\mathcal{T}(\cdot,\Theta_0,\Lambda_0))$ on $[0,\tau_f]$.

Let a time $\tilde\tau$ be chosen from interval $[\bar\tau,\tau_f)$, and let $B$ represent the system state at $\tilde\tau$. Define $\gamma$ as the segment of the extremal trajectory $\Pi(\mathcal{T}^\mu(\cdot,\Theta_0,\Lambda_0))$ from $B$ to the target. Similarly, let $\hat\gamma$ denote the smooth concatenation of the extremal trajectory $\Pi(\mathcal{T}^\mu(\cdot,\Theta_0,\Lambda_0))$ from $B$ and $A$ and the trajectory of the dashed curve. Clearly, the cost consumed by the pursuer when following $\gamma$ is identical to the cost associated with traversing $\hat\gamma$. However, along the modified trajectory $\hat\gamma$, a discontinuity in the control input arises at point $A$. This contradicts with the necessary condition in Eq.~(\ref{eq:u}) in which the control is continuous. Thus, there is another trajectory in the neighborhood of $\hat{\gamma}$ from $B$ to the target so that the cost is smaller, indicating that there is a trajectory requiring smaller control effort than the extremal trajectory $\gamma$. This contradicts with the assumption, completing the proof. 
$\hfill\square$

\subsection{Proof of Lemma \ref{lemma:rotation}}\label{app:C}
\setcounter{equation}{0}
\renewcommand{\theequation}{C\arabic{equation}}
Given any feasible state $(r_c,\lambda_c,\theta_{Pc})$ with a heading angle of the target $\theta_T\in(0,2\pi)$, and a speed ratio $\mu\in(0,1)$, there exists an optimal trajectory $(r(t),\lambda(t),\theta_P(t))$ for $t\in[0,t_f]$ so that $(r_c,\lambda_c,\theta_{Pc})=(r(0),\lambda(0),\theta_P(0))$. Notice that we have
\begin{equation}
	\label{eq:ode2}
	\begin{cases}
		\displaystyle{\frac{d}{dt}}r(t)=-\cos\sigma_P(t)+\mu\cos\sigma_T(t)\\
		\displaystyle{\frac{d}{dt}}\lambda(t)=\displaystyle{\frac{-\sin\sigma_P(t)+\mu\sin\sigma_T(t)}{r(t)}}\\
		\displaystyle{\frac{d}{dt}}\theta_P(t)=u(t)\\
		%&\dot \theta_T=0
	\end{cases}
\end{equation}
Set $\tau=t$, and let
\begin{equation}\label{eq:C2}
	\begin{cases}
		\overline{r}(\tau)=r(\tau)\\
		\overline\lambda(\tau) = \lambda(\tau)+\theta_T\\
		\overline\theta_P(\tau)=\theta_P(\tau)+\theta_T\\		
	\end{cases}
\end{equation}
Then, according to Eq.~(\ref{eq:ode2}) and Eq.~(\ref{eq:C2}), we have
\begin{equation*}
	%\label{equation:ode1}
	\begin{cases}
		\displaystyle{\frac{d}{d\tau}}\overline{r}(\tau)=-\cos(\overline\theta_P(\tau)-\overline\lambda(\tau))+\mu\cos(0-\overline\lambda(\tau))\\
		\displaystyle{\frac{d}{d\tau}}\overline\sigma_T(\tau)=\displaystyle{\frac{-\sin(\overline\theta_P(\tau)-\overline\lambda(\tau))+\mu\sin(0-\overline\lambda(\tau))}{\overline{r}(\tau)}}\\
		\displaystyle{\frac{d}{d\tau}}\overline\theta_P(\tau)=\overline{u}(\tau)\\
	\end{cases}
\end{equation*}
Since $u(t)$ represents the optimal control corresponding to the trajectory $(r(t),\lambda(t),\theta_P(t))$ for $t\in[0,t_f]$ with the heading angle of the target $\theta_T$, it follows that $\overline{u}(\tau)$ serves as the optimal control for the trajectory $(\overline{r}(\tau),\overline\lambda(\tau),\overline\theta_P(\tau))$ for $\tau\in[0,t_f]$ with a heading angle of the target $\theta_T=0$. It is apparent that a constant increment on the LOS angle $\lambda$ and pursuer's heading angle $\theta_P$ will not change the trajectory. Thus, we have
\begin{equation*}
	u^*(r(t),\lambda(t),\theta_P(t),\theta_T,\mu)=u^*(\overline{r}(\tau),\overline\lambda(\tau),\overline\theta_P(\tau),0,\mu)
\end{equation*}
This completes the proof.$\hfill\square$

\subsection{Proof of Lemma \ref{lemma5}}\label{app:D}
\setcounter{equation}{0}
\setcounter{figure}{0}
\renewcommand{\theequation}{D\arabic{equation}}
\renewcommand{\thefigure}{D\arabic{figure}}
Given any feasible state $(r_c,\lambda_c,\theta_{Pc})$, a heading angle of the target $\theta_T=0$, and any speed ratio $\mu\in(0,1)$, there exists an optimal trajectory $(r(t),\lambda(t),\theta_P(t))$ for $t\in[0,t_f]$ so that $(r_c,\lambda_c,\theta_{Pc})=(r(0),\lambda(0),\theta_P(0))$. Notice that we have
\begin{equation}
	\label{eq:ode3}
	\begin{cases}
		\displaystyle{\frac{d}{dt}}r(t)=-\cos\sigma_P(t)+\mu\cos\sigma_T(t)\\
		\displaystyle{\frac{d}{dt}}\lambda(t)=\displaystyle{\frac{-\sin\sigma_P(t)+\mu\sin\sigma_T(t)}{r(t)}}\\
		\displaystyle{\frac{d}{dt}}\theta_P(t)=u(t)\\
		%&\dot \theta_T=0
	\end{cases}
\end{equation}
Set a positive scalar $s>0$ and $\tau=t/s$, and let
\begin{equation}\label{eq:1/s}
	\begin{cases}
		\overline{r}(\tau)=\displaystyle{\frac{1}{s}}r(s\tau)\\
		\overline\lambda(\tau) = \lambda(s\tau)\\
		\overline\theta_P(\tau)=\theta_P(s\tau)\\		
	\end{cases}
\end{equation}
Then, according to Eq.~(\ref{eq:ode3}) and Eq.~(\ref{eq:1/s}), we have
\begin{equation*}
	%\label{equation:ode1}
	\begin{cases}
		\displaystyle{\frac{d}{d\tau}}\overline{r}(\tau)=-\cos(\overline\theta_P(\tau)-\overline\lambda(\tau))+\mu\cos(\theta_T-\overline\lambda(\tau))\\
		\displaystyle{\frac{d}{d\tau}}\overline\sigma_T(\tau)=\displaystyle{\frac{-\sin(\overline\theta_P(\tau)-\overline\lambda(\tau))+\mu\sin(\theta_T-\overline\lambda(\tau))}{\overline{r}(\tau)}}\\
		\displaystyle{\frac{d}{d\tau}}\overline\theta_P(\tau)=su(t)\\
	\end{cases}
\end{equation*}
It follows that the trajectory $(\overline{r}(\tau),\overline\lambda(\tau),\overline\theta_P(\tau))$ for $\tau\in[0,t_f/s]$ is also the optimal trajectory of the OCP. According to Lemma \ref{lemma:ex_tra}, there exists an extremal trajectory $\Pi(\mathcal{T}^\mu(\tau,\Theta_0,\Lambda_0))$ for $\tau\in[0,t_f/s]$ so that
\begin{equation*}
	\begin{cases}
		\overline{r}(\tau)=R^\mu(\tau,\Theta_0,\Lambda_0)\\
		\overline\lambda(\tau)=\Lambda^\mu(\tau,\Theta_0,\Lambda_0)\\	
		\overline\theta_P(\tau)=\Theta^\mu(\tau,\Theta_0,\Lambda_0)\\
	\end{cases}
\end{equation*}
%Thus, we have $\Pi(\mathcal{T}^\mu(\tau,\Theta_0,\Lambda_0))\in\mathcal{F}^\mu$, which completes the proof of the first statement.

Since $u(t)$ represents the optimal control corresponding to the trajectory $( r(t),\lambda(t),\theta_P(t))$ for $t\in[0,t_f]$, it follows that $su(s\tau)$ serves as the optimal control for the trajectory $(\overline{r}(\tau),\overline\lambda(\tau),\overline\theta_P(\tau))$ for $\tau\in[0,t_f/s]$. Thus, we have
\begin{equation*}
	\begin{split}
		&C(\overline{r}(\tau),\overline\lambda(\tau),\overline\theta_P(\tau),\mu)=sC({r}(t),\lambda(t),\theta_P(t),\mu)
	\end{split}
\end{equation*}
which further indicates
\begin{equation*}
	\begin{split}
		C({r}(t),\lambda(t),\theta_P(t),\mu)=\frac{1}{s}C(\overline{r}(\tau),\overline\lambda(\tau),\overline\theta_P(\tau),\mu)
	\end{split}
\end{equation*}
This completes the proof.$\hfill\square$

\section*{Acknowledgments}
This research was supported by the National Natural Science Foundation of China under grant Nos. 61903331 and 62088101.

\bibliography{sample}

%\begin{IEEEbiography}[{\includegraphics[width=1in,height=1.25in,clip,keepaspectratio]{Me.png}}]{Han Wang}{\space}received his B.S. degrees in the major of Detection, Guidance and Control from Xidian University, People's Republic of China, in 2020. He has been a Ph.D. candidate in the major of Aerospace Science and Technology, Zhejiang University since 2020. His research interest is in the development of machine learning techniques for real-time generation of optimal guidance command in aerospace engineering.
%\end{IEEEbiography}%
	
%\begin{IEEEbiography}[{\includegraphics[width=1in,height=1.25in,clip,keepaspectratio]{cz.jpg}}]{Zheng Chen}{\space}received the Ph.D. degree in Applied Mathematics from the University Paris-Saclay in 2016, the M.S. and B.S. in Aerospace Engineering from Northwestern Polytechnical University in 2013 and 2010, respectively. He is currently a professor with the School of Aeronautics and Astronautics at Zhejiang University. His current research interests revolve around flight mechanics and nonlinear optimal guidance  in aerospace engineering.
%\end{IEEEbiography}

\end{document}